\documentclass[a4paper,11pt]{amsart}

\usepackage{fancyhdr}
\usepackage{amsmath,inputenc,euscript,amssymb}
\usepackage{fullpage}
\usepackage{bbm}
\usepackage{color}
\usepackage[colorlinks=true, linkcolor=magenta, citecolor=cyan, urlcolor=blue]{hyperref}
\usepackage{datetime}

\font\bb=msbm10 at 11pt
\font\bbi=msbm10 at 8pt
\def\N{\hbox{\bb N}}

\def\R{\hbox{\bb R}}

\def\Ri{\hbox{\bbi R}}

\newtheorem{theo}{Theorem}[section]
\newtheorem{pro}[theo]{Proposition}

\newtheorem{cor}[theo]{Corollary}

\newtheorem{rem}[theo]{Remark}

\title{Upper Bounds for the   Blow-up Time of the 2-D  Parabolic-Elliptic Patlack-Keller-Segel  Model of Chemotaxis}

 \subjclass[2010]{
35Q92, 35B44, 35K55, 35M31, 92C17.
}

\keywords{Patlack-Keller-Segel system; chemotaxis; critical mass; blow-up time upper bound.}

\author{Patrick Maheux}

\address{Institut Denis Poisson (IDP) \\ 
Universit\'e d'Orl\'eans, Universit\'e de Tours, CNRS -- UMR 7013 \\
B\^atiment de Math\'ematiques, B.P. 6759\\
45067~Orl\'eans Cedex 2, France }
\email{Patrick.Maheux@univ-orleans.fr}

 \thanks{ 
This work was carried out while the author benefited from a delegation at the CNRS 
and a CRCT during the first semester of 2021 and the first semester of 2023, respectively. 
}
 
\begin{document}

\maketitle

\begin{center}
\today  
\end{center}

\begin{abstract}
In this paper, we obtain  upper bounds for  the critical time $T^*$  of the  blow-up for the  parabolic-elliptic Patlak-Keller-Segel 
system on  the  2D-Euclidean space. 
No  moment condition or/and entropy condition are required on the initial data; only the usual assumptions of non-negativity and 
finiteness of the  mass are assumed.
The result is expressed not only in terms of the supercritical  mass $M> 8\pi$, but also in terms of  the {\it shape} of the initial data.
\end{abstract}

\footnote{I would like to warmly thank my colleague Vittoria Pierfelice for introducing me to the subject of Patlak-Keller-Segel  Models for Chemotaxis \cite{MP1}.}

\tableofcontents

\section{General upper bound on $T^*$ for the (PKS) system}\label{intro}

In this paper, we consider the following  Patlak-Keller-Segel (PKS) model for chemotaxis in the whole plane ${\R}^2$ \cite{KellerSegel,Patlak,Jager}:
\begin{equation}\label{pks1}
\frac{\partial}{\partial t}n(x,t)=\Delta n(x,t)-\nabla.(n\nabla c)(x,t), \quad x\in \R^2,  \; t>0, 
\end{equation}
\begin{equation}\label{pks0}
(-\Delta)c(x,t)=n(x,t),  \quad x\in \R^2,  \; t>0, 
\end{equation}
 \begin{equation}\label{init}
 n(0,x)=n_0(x),  \quad x\in \R^2,  \quad 0\leq n_0\in L^1(\R^2),
\end{equation}
 where $n_0$ denotes the initial data.
The function $n(x,t)$ represents the density of  cells, and $c(x,t)$ 
 the density of the chemical secreted by themselves that attracts them and makes them to aggregate.	

There is a huge bibliography on  this type of systems and  its  variants.  
The geometric framework is usually  $\R^n$ or bounded domains on $\R^n$ with different kind of results depending on the dimension   $n=1,2$ or $n\geq 3$
and assumptions of the initial data. 
Note that  there are also  some  recent works for specific manifolds with curvature, see  \cite{MP1,Nagy,PTNN} for instance. 
Diverse analytic settings for the existence of the solutions for the (PKS) systems have been introduced as classical solutions, strong solutions, entropy solutions,  
various notions of weak solutions; see  \cite{Bil} and references therein. 
Stochastic Keller-Segel particle models have been also introduced for the existence and the uniqueness of solution of (PKS) systems,  
see   for instance \cite{CaPe,Tomc, FJ} and  \cite{Fat} for numerical simulations.  

For various versions of the (PKS) systems,  blow-up time bounds have been obtained more or less explicitly.
For instance, upper bounds on the maximal existence time $T^*$ have been mentioned on  the disk  $\Omega=B_R$ ($R>0$) of  $\R^2$  in \cite[Thm 1.1, eq. (2.1)]{CiSt1}
for the parabolic-parabolic (PKS) systems (but not explicitly given).
Always for the parabolic-parabolic (PKS) systems, we can also mention  \cite[Thm. 1.2]{CC1} for a study  on the whole plane under a second moment condition for $n_0$.
In higher dimensions, we can see for instance  \cite[Section 4.1]{CC3}, also \cite[Thm 1.1, Lemma 4.1]{CiSt2}.
In many papers, the bounds on $T^*$ are not
explicit and/or  impose some extra-assumptions on the initial data.
 It is not our intention to cite  all the possible articles on this subject. So,
 we shall limit ourselves to mention a few of the papers on (PKS) systems related 
to bounds on the maximal existence time $T^*$, especially on $\R^2$ which is the  setting of study in  this paper.
We apologize in advance for not mentioning more articles, in particular those  dealing   with  generalizations of  
(PKS) systems and/or  studying  other types of questions for these systems. More references can be found in  the recent book by P.Biler \cite{Bil}.

Among many papers on the subject, we  can mention some of them historically at the beginning of the study of the  blowup in finite time, on
chemotactic collapse (convergence to Dirac measure plus possibly a $L^1$ function at the critical time $T^*$), and  on the profile of solutions  at $t=T^*$ 
on bounded domain, or on the whole $\R^2$ for (radial) solutions of (PKS) systems: see
 \cite{HV0,HV1,HV2, Nagai1, Nagai2, Nagai3}
 and, more recently, \cite{Mizo,CGMN} and references therein. 
A variant of the (PKS) system on  a ball of $\R^n$, $n\geq 2$ is studied in \cite{Winkler} for which  a novel type of critical mass phenomenon 
 for radially symmetric initial data and  linked to the concept of concentration comparison  is described implying  that the solutions blow up in finite time. 
 Note that the concept of concentration comparison takes into account in some way of the shape of the initial data $n_0$, see  Theorem 1.1 in  \cite{Winkler}.
 See also \cite{BCKZ} for a model of (PKS) system with consumption term on $\R^2$.
We shall not involved in such refinement but our estimates depend also on the shape of the initial data $n_0$.

It seems that no  specific attempt has been  made to  estimate the critical time $T^*$ over the whole space $\R^2$ 
with minimal assumptions on the initial data.
In this paper, we make some progress in this direction by  providing  information on upper bounds on $T^*$.
But we shall not be concerned with lower bounds on $T^*$ with this generality in this article. 
Nevertheless, we  notice  that in \cite[p.355]{koz-sug}, and for some specific situations, some  lower bounds on $T^*$ have been obtained over the whole space $\R^2$.
On bounded convex domains  of $\R^2$ and $\R^3$, lower bounds on the  maximal existence time $T^*$ in case of blow-up 
(with  various notions of blow-up) have been investigated first by Payne-Song,  see \cite{PaSo1,PaSo2}.
 See also \cite[Th. 3.4]{MNV} for a recent work dealing with  more general (PKS) systems on bounded smooth domains of $\R^n$ with $n\geq 1$.
In  \cite{FMV}, some numerical simulations are given  and compared with analytic  bounds of \cite{PaSo1,PaSo2}.
The study  of lower bounds on maximal existence  time $T^*$ has  regained some interest more recently, 
in particular for  some variants of the parabolic-parabolic version of (PKS) on bounded domains with smooth boundary;
see \cite{TaPi,TeTo,MNV}  and references therein. 

\vskip0.3cm

Let us return to the main  subject of this paper.
In the  framework of weak solutions, it is proved in  \cite[Cor.2.2]{BDP}, 
 under finite second moment condition and finite entropy of $n_0$, the following result.
For $M>8\pi$,
the maximal existence time $T ^*$ of the solutions
 is bounded as follows,
\begin{equation}\label{bupbdp}
 T^*\leq \frac{2\pi I(0)}{M(M-8\pi)},
\end{equation}
 where $I(0)=\int_{\Ri^2} \vert x\vert^2 n_0(x)\, dx$ is the second moment of $n_0$.
 The inequality \eqref{bupbdp} expresses  a blow-up of the solutions of the (PKS) system, i.e. $T^*<\infty$. 
 Note that this result also holds true with $I(0)$ replaced by
 $$
 MV(0)=\int_{\Ri^2} \vert x-B_0\vert^2 n_0(x)\, dx,
 $$ where 
 $V(0)$ is the  variance of the initial data $n_0$ defined by 
 $$
 V(0)=V_2(n_0)=\frac{1}{M}\int_{\Ri^2} \vert x-B_0\vert^2 n_0(x)\, dx,
 $$
and $B_0$  is given by 
 $$
 B_0=B_0(t)=\frac{1}{M} \int_{\Ri^2}   x.n_t(x)\, dx=\frac{1}{M}\int_{\Ri^2}   x. n_0(x)\, dx \in \R^2, \quad t\in (0,T^*).
 $$ 
 Here, the measure  $\frac{n_0(x)}{M}\, dx$ can be seen as a probability measure since $n_0$ is non-negative and non identically zero.
 For a study of (PKS) system with Borel measures as initial data, see \cite{Bed-Mas}.
 It can be easily seen that the barycenter  $B_0(t)$ is (formally) independent of $t\in (0,T^*)$ for any solution $(n_t)$  of the (PKS) system.
 Thus we have the following  easy improvement of \eqref{bupbdp}, namely
\begin{equation}\label{bupbdp2}
 T^*\leq T^*_v:= \frac{2\pi V(0)}{M-8\pi}.
\end{equation}
Let us recall the formal  arguments of the proof  for these inequalities \eqref{bupbdp} and \eqref{bupbdp2}. It is first enough to show  \eqref{bupbdp2}.
We have the following  exact formula,
 $$
 I^{\prime}(s):=\frac{d}{ds}\left( \int_{\Ri^2} \vert x\vert^2 n_s(x)\, dx\right)
 =4M\left(1-\frac{M}{8\pi}\right), \quad s\in (0,T^*),
 $$
 for  any weak solution  $(n_t)$   of (PKS) system with  mass  $M=\vert\vert n_0\vert\vert_1$ such that 
 $$ \int_{\Ri^2} \vert x\vert^2 n_0(x)\, dx< +\infty,
 $$
 see \cite[Lemma 2.1]{BDP}.
 Let  $V(s)=\frac{1}{M}\int_{\Ri^2} \vert x-B_0\vert^2 n_s(x)\, dx
 =
 \frac{I(s)}{M}-\vert B_0\vert^2$.
 The derivative  just above   can also  be written as 
 $$
 V^{\prime}(s)=\frac{1}{M}  I^{\prime}(s)
 =4\left(1-\frac{M}{8\pi}\right), \quad s\in (0,T^*).
 $$
By  integration with respect to $s$, this leads   to the next formula 
 $$
 V(t)=V(0)+4\left(1-\frac{M}{8\pi}\right)t, \quad 0<t<T^*.
 $$
 By maximal principle $n_t\geq 0$, so  $V(t)$  is non-negative for  all $0<t<T^*$. Thus, we obtain  
 $$
  V(0)-4t\left(\frac{M-8\pi }{8\pi}\right)\geq 0,  \quad 0<t<T^*.
  $$
  Hence, we conclude that
  $ t\left(\frac{M-8\pi }{2\pi}\right)\leq V(0)$ for all  $ 0<t<T^*$. Letting $t$ goes to $T^*$, we deduce \eqref{bupbdp2}
when $M>8\pi$.
Now, note that \eqref{bupbdp2} implies  \eqref{bupbdp} because
$$
MV(0)=I(0)- M \vert B_0\vert^2\leq I(0).
$$
Note that usually the barycenter is rarely mentioned in most papers due probably to the fact that the initial data 
is often  assumed  to be radially symmetric, hence $B_0=0$. 
But we  see some advantages to consider  barycenters $B_0\neq 0$ and the variance in the formulation of some of our results.
Recall that the variance   measures the spread of the data around the mean. 
In particular,  when  we mention the variance of the density $\frac{n_0(x)}{M}\, dx$.
Later on, this result has been generalized  in \cite[Th.2]{koz-sug} 
where  the second equation \eqref{pks0} of the (PKS) system is replaced by 
$$
(-\Delta)c(x,t)+\gamma c(x,t)=n(x,t),  \quad x\in \R^2,  \; t>0, 
$$
with $\gamma\geq 0$. When $\gamma >0$, we say that we have a {\it consumption}  term.
In that case,  an additional assumption  on $I(0)$ of the form  $I(0)\leq h_1(M)$ for some function $h_1$ of the mass $M$ is imposed to obtain the blow-up
in the supercritical case.
 
A result  has recently  been obtained  for the analogue of the (PKS) system on the 2-D  hyperbolic space. 
More precisely, with an appropriate definition of the moment  $I(0)$ for the initial data $n_0$, and  under  a similar additional condition 
of the form $I(0)\leq h(M)$  for some function $h$ of the mass $M$, 
a blow-up is proved for the solution  of (PKS) system.
Moreover, the maximal time $T^*$ is bounded by an explicit function of the mass $M$  and $I(0)$; analogue to \eqref{bupbdp2},  
see \cite{MP1}.
This last situation shares some similarity with the Euclidean case with consumption, (i.e. $\gamma>0$), 
certainly due to the spectral gap of the Laplacian on the  hyperbolic space.
Recently, local criteria have been  introduced for blow-up of radial solutions in 2-dimensional chemotaxis models have been obtained by P. Biler et al.,
using weighted averages on disks.
See \cite{BKZ}, where the Laplacian $\Delta$ is replaced by   some powers of the Laplacian $(-\Delta)^{\alpha/2}$.
See also  \cite{BCKZ}, when  the consumption term with $\gamma >0$ is considered.

Let us mention also the  case of critical mass $M=8\pi$ on  the plane  for which the critical time is $T^*=+\infty$, i.e. solutions are global in time,  
when assumptions of finite 2-moment and entropy are made on the initial condition $n_0$. 
Moreover, the solutions blow up as the  delta Dirac  $8\pi\delta_{z_0}$ at the center of mass $z_0=MB_0$ when $t\rightarrow T^*=+\infty$, 
see \cite{BCM} [Th.1.3, p.1453],  see also \cite{HV2} and the recent long preprint \cite{DDPDMW}. 
In our study, we shall see that the upper bound $T_c^*$ of $T^*$ is consistent with this result when $M\rightarrow 8\pi^+$. 
See Remark \ref{rmrm} below.
\vskip0.3cm

The aim of the present work is to provide an upper bound on the maximal existence time $T^*$ in the case of blow-up $M>8\pi$
by removing the second moment condition and/or  the entropy condition on  the initial  data $n_0$. We only assume that $n_0$ is non-negative and integrable on $\R^2$.

Throughout this paper, we shall use   the definition of mild solution $n(t)$ of the (PKS) system on $\R^2$ taken from
 \cite[p.392]{Dongyi Wei}.
Let $T>0$ be fixed.
We say that   $n(t)=n(.,t)$, $0<t<T$,  is a  mild solution of the (PKS) system if,
$$
n\in C_{w}([0,T), L^1(\R^2)),\quad \sup_{t\in (0,T)}
t^{\frac{1}{4}}
\vert\vert n(t)\vert\vert_{L^{\frac{4}{3}}}
< +\infty,
$$
and $n(t)$ satisfies the following Duhamel integration equation for all $t\in (0,T)$,
\begin{equation}\label{duhamel}
n(t)=e^{t\Delta}n_0-\int_0^t e^{(t-s)\Delta}div(n(s)\nabla c(s))\,ds,
\end{equation}
with $-\Delta c(t)=n(t)$ in the sense
$$
\nabla c(x,t)=-\frac{1}{2\pi} \int_{\Ri^2}\frac{x-y}{\vert x-y\vert^2} n(t,y)\,dy,
 \quad x\in \R^2,  \; t\in (0,T).
$$
Here,  the space  $C_{w}([0,T), L^1(\R^2))$ is the space of (weakly) continuous functions with values in $L^1(\R^2)$.
Recall that, for any   initial data $0\leq n_0\in L^1(\R^2)$, there exists  a  
mild solution $(n(t))_{0<t<T^*}$ of the (PKS) system where 
$T^*\in (0,+\infty]$ denotes the maximal  existence time  of the solution (see \cite{Dongyi Wei,Bil}).
For the case of  existence of (global)  weak solutions with measures as initial data on $\R^2$, we can refer to  the recent paper \cite{FT}
where a simple proof of this fact is  given.

Recently,   Dongyi Wei   has  proved  the following remarkable dichotomy result.

\begin{theo} \cite[Th. 1.1 p.390]{Dongyi Wei}\label{dongyiweithm}
Assume that $0\leq n_0\in L^1(\R^2)$ and $M=\vert\vert n_0\vert\vert_1$.
Let $T^*$ be the maximal  existence time  of the mild solution.
Then $T^*=+\infty$ if and only if $M\leq 8\pi$.
\end{theo}

As already described above, an additional  finite second moment condition on the initial data $n_0$  is  often used
to imply by a simple virial argument the blow-up of the solution and to provide  an explicit upper bound on the maximal existence time $T^*$ of the solution 
 when $M>8\pi$, see  for instance \cite{BDP, koz-sug}.
The main novelty of the result  of D. Wei is that no moment condition (nor entropy)  is imposed  on the initial data $n_0$
other than the  minimal ones \eqref{init}.
In other words, his  result can be rephrased as follows:  the (mild) solution blows up, i.e. $T^*<+\infty$,  if and only $M>8\pi$. 
But it  seems that there are  no explicit  estimates  of $T^*$  deduced from the study carried out in  \cite{Dongyi Wei}. 
\vskip0.3cm 

Before describing the main results  of this paper, we make some general comments 
on the approach  which consists in replacing the unbounded weight $\vert x \vert^2$ 
appearing in the definition of the second moment by a family of bounded weights: 
Gaussian  weights (i.e. heat kernels).
We start our discussion by some simple facts related to the heat equation which is known to play
an important role in the study of (PKS) systems, in particular by the Duhamel formula \eqref{duhamel}. 
This importance   is due to the fact that the heat kernel  is related to  the linear part of  the (PKS) equation 
\eqref{pks1}, i.e. to the Laplacian $\Delta$.
Let us denote by  $p_s$  be the heat kernel on $\R^2$ defined by 
$$
p_s(x)=({4\pi s})^{-1} \exp\left(-\frac{\vert x\vert^2}{4s}\right)
$$
for all  $x\in \R^2$ and  all $s>0$.
We  also denote by $H_z(s)=H_{z,n_0}(s)$  the following function 
$$
H_{z,n_0}(s):=4\pi s \, e^{s\Delta}n_0(z)=
4\pi s \,p_s\star n_0(z)=
\int_{\Ri^2} \exp\left(-\frac{\vert x -z\vert^2}{4s}\right)n_0(x)\, dx
$$
for all   $z\in \R^2$ and  all $s>0$.
The family of functions $H_{z,n_0}$ can be seen as particular weighted integrals.
The general  interest of studying  weighted  integrals of the form
$$
I_{w}(t):=\int_{\Ri^2}  w(x)n_t(x)\, dx, \quad 0<t<T^*, 
$$
for some  positive weight $w$ is to {\it detect} and {\it measure in average}
the effect of the evolution  in time of  the density $(n_t)_{ 0<t<T^*}$ of cells.
If  the chosen  weight  $w$ is constant, the only  fact that  can be observed is that the mass of  cells is preserved.
Indeed, the mass  of $n(x,t)$ given by $I_{w}(t)$ with $w\equiv 1$, i.e.,
$$
M(t)=\int_{\Ri^2} n(x,t)\, dx=M(0)= M=\int_{\Ri^2} n_0(x)\, dx,
$$
 is preserved along the time evolution of the (PKS)  system.
When one  considers  the weight $w(x)=\vert x\vert^2$, such a  weight imposes 
an additional condition on $n_0$, i.e.,  the finiteness of the second-moment
$\int \vert x\vert^2n_0(x)\, dx$, which leads automatically   to a lack of generality. 
See \eqref{bupbdp} above  for the corresponding bound on  $T^*$.

Other choices of $w$  may   detect other  behaviours of the solution $(n_t)$ in the sense that 
 the quantity  $I_{w}(t)$ may evolve with $t$  in  diverse ways, see for instance \cite{Bil} eq. 5.3.13 p.160. 
 This is what happens when we choose the heat kernel as weight  to describe a two  parameters  family of weights, 
 namely $w_{s,z}(x)=  {4\pi }sp_s(x-z)$ with  $s>0$ and  $z\in \R^2$ be fixed.
The use of  the heat kernel as weight  has at least three  features.
The first interest is that the  weight 
$$
w_{s,z}(x)=  {4\pi } sp_s(x-z)
= 
\exp\left(-\frac{\vert x -z\vert^2}{4s}\right)
$$
is nowhere constant, smooth  and uniformly bounded by $1$  for all  $s>0$ and $z\in \R^2$. 
In particular, the function $(s,z)\mapsto H_{z, n_0}(s)$ defined above is always  finite, more precisely  
$H_z(s)\leq M$ where $M$ is  the mass  of $n_0$ which  is assumed to be finite. 
In this way, no additional condition on $n_0$ than $n_0\in L^1(\R^2)$ is necessary to study the quantity $H_{z, n_0}(s)$.
The  second argument is that   the  function  $s\mapsto H_{z, n_0}(s)$ is   continuous, 
monotonically increasing with $s$, and  its range is  $(0,M)$ for a non-negative non-zero initial data $n_0$.
 Furthermore, choosing the point  $z \in \R^2$ allows us to emphasize the behavior of the initial data $n_0$  around this point
due to the rapid decay of the heat kernel at infinity.
Finally, the third reason and probably the most important  is that the heat  kernel has a strong relationship with the (PKS) system. Indeed, 
the simplest linearized form of the non-linear equation \eqref{pks1}, i.e. without the aggregation term $-\nabla.(n\nabla c)$,
is the heat equation, i.e., $\frac{\partial}{\partial t}u(x,t)=\Delta u(x,t)$.
\vskip0.3cm

The aim of this paper is to provide  a general  upper  bound on the maximal  existence time   $T^*$ of  mild solutions of (PKS) system 
in the supercritical  case  $M>8\pi$ without further conditions on the initial data $n_0$ other  than $0\leq n_0\in L^1(\R^2)$. 
This is achieved  by Theorem \ref{tcalphazero}  below, which proposes a general   formula for an upper bound $T_c^*(n_0)$ of  $T^*$.
Later on, in Theorem \ref{tclaplace},  we specify the form of the bound on $T_c^*(n_0)$ for  radially symmetric initial data, which allows us 
to use the Laplace transform.
The main interest of the theorem just below  is to evaluate an upper bound on the maximal existence time $T^*$.
No additional conditions, as a moment condition and/or an entropy condition, are  imposed on the initial data $n_0$.

These estimates could be used as time bounds for numerical simulations on the time interval $[0,T)$ 
 observing with certainty the blow-up  phenomenon on the plane by choosing $T \geq T_c^*(n_0)\geq T^*$, see Theorem \ref{tcnormal}.
But for  a reasonable simulation,  i.e on a bounded domain $\Omega$ with smooth boundary $\partial \Omega$,  it would be necessary to consider the case where the cells are located at time $t = 0$  in a  set at a very  large distance from the  boundary of  the domain  to avoid  possible  boundary  effects.
Recall that $n_0$ represents the repartition of the cells at time $t=0$. 
So, the condition will certainly be described   by a  control of the form
${\rm dist}({\rm supp} \, n_0, \partial \Omega)\geq C$ with $C>0$ large enough.
\vskip0,3cm

Here is the first  main result of this paper.

\begin{theo}\label{tcalphazero}
Let  $0\leq n_0\in L^1(\R^2)$  of  mass  $M=\int_{\Ri^2} n_0(x)\, dx>8\pi$ and $T^*$  be the maximal   existence  time 
of a mild   solution $(n_t)$ of  the Patlack-Keller-Segel system \eqref{pks1}-\eqref{pks0}-\eqref{init}.
 Then the following statements hold true.
\begin{enumerate}
\item
The maximal existence time $T^*$  is finite and satisfies the following estimate,
\begin{equation}\label{heatcondit}
 H_{z,n_0}(T^*):=\int_{\Ri^2} \exp\left(-\frac{\vert x -z\vert^2}{4T^*}\right)n_0(x)\, dx\leq
\frac{2M^2}{3M-8\pi},
\end{equation}
for all  $z\in \R^2$.
In particular,  $T^*$  is finite.
\vskip0.3cm

Or equivalently,  for all $0<s<T^*$
 and all $z\in \R^2$,
\begin{equation}\label{heatcondit2}
 H_{z,n_0}(s):= \int_{\Ri^2} \exp\left(-\frac{\vert x -z\vert^2}{4s}\right)n_0(x)\, dx\leq
\frac{2M^2}{3M-8\pi}.
\end{equation}
\item
For $z\in \R^2$, let  $T_{c,z}^*(n_0):=H_{z,n_0}^{-1}\left(\frac{2M^2}{3M-8\pi}\right)$. Then we have
\begin{equation}\label{tcnormal}
T^*\leq 
T_c^*(n_0):=
\inf_{z\in \Ri^2}T_{c,z}^*(n_0)<+\infty.
\end{equation}
\item
Moreover, if there exists $z_0\in\R^2$ such that 
$H_{z,n_0}(s)\leq H_{z_0,n_0}(s)$ for all $s>0$ and all $z\in \R^2$, then $T_c^*(n_0)=
 T_{c,z_0}^*(n_0)$. In particular, if $n_0$ is a non-increasing radially symmetric  integrable  function, then  
 $$
 T_c^*(n_0)=
 T_{c,{z=0}}^*(n_0)=H_{0,n_0}^{-1}\left(\frac{2M^2}{3M-8\pi}\right).
 $$
\item
The value $T_c^*(n_0)$ is translation-invariant in the sense that
$T_c^*(n_0)=T_c^*(m_0)$ for any $m_0$ of the form $m_0(x)=n_0(x+z_0)$, $x\in \R^2$ for fixed $z_0\in \R^2$.
\end{enumerate}
\end{theo}

Before giving the proof of the theorem just above, we make  several comments.  

\begin{rem}
 At first sight, the inequality \eqref{heatcondit2} comes as a surprise since the full range of $H_z(s):=H_{z,n_0}(s)$, $s>0$,  is a priori  $(0,M)$
 whenever $n_0\neq 0$ and $n_0\in L^1(\R^2)$. Note that $H_z(s)$ only depends on the Gaussian and any initial data  $n_0$.
 The inequality \eqref{heatcondit2}  is a limitation of this range since we easily see that
 $$
 \frac{2M^2}{3M-8\pi}
 <M,
 $$ 
 when $M>8\pi$. Hence, we  deduce qualitatively  that $T^*$ is finite. 
 Indeed, if $T^*=+\infty$  we can take the limit as $s$ goes to infinity in \eqref{heatcondit2}, 
 and by monotone convergence we get 
 $$
M= \lim_{s\rightarrow +\infty} H_{z,n_0}(s)\leq \frac{2M^2}{3M-8\pi},
 $$
 for any $z\in \R^2$.
Contradiction.
 This argument also  implies  that all the quantities  $T_{c,z}^*(n_0)$  defined in Theorem \ref{tcalphazero}
 are finite  since $s\rightarrow H_{z,n_0}(s)$ is strictly increasing, hence it is invertible on its range $(0,M)$
 which contains $ \frac{2M^2}{3M-8\pi}$.
 In particular, the inequality $T^* \leq T_{c,z}^*(n_0)$  holds true  for all $z\in \R^2$ when $M>8\pi$ which allows us  to optimize over $z$.
 Thus, the lower bound of the upper bounds   $T_{c,z}^*(n_0)$ over $z$ will provide  an  upper bound on $T^*$.
\end{rem}

\begin{rem}\label{rmrm}
The constant $\frac{2M^2}{3M-8\pi}$ in the inequality  \eqref{heatcondit}-\eqref{heatcondit2}  is sharp  for the supercritical case  in the sense that  the inequality 
$R(M):=\frac{2M^2}{3M-8\pi}<M$  holds for all  $M>8\pi$, and this inequality becomes an equality when $M\rightarrow 8\pi^+$. 
Note that this sharpness  is essential for the proof of  the dichotomy  result of Theorem 1.1 in
 \cite{Dongyi Wei}, and  consequently also for  Theorem \ref{tcalphazero} of this article. Note that when $M=8\pi$, then $R(M)=M=8\pi$, and formally 
 $$
 T^*= T_{c}^*(n_0)= T_{c,z}^*(n_0)=\lim_{M\rightarrow 8\pi^+} H_{z,n_0}^{-1}\left(\frac{2M^2}{3M-8\pi}\right)
 =
 +\infty
 $$
 for all $z\in \R^2$. Recall that when $M=8\pi$ then we have $ T^*= +\infty$,  see  \cite[Th. 1.1 p.390]{Dongyi Wei},
recalled in  Theorem \ref {dongyiweithm} of this paper.
 This  is consistent with a former  existence result  of  \cite{BCM} of global solutions for the critical mass $M=8\pi$
 when  2-moment  and  entropy of the initial data $n_0$ are assumed to be finite.
Moreover, in  \cite{BCM}, a more precise study  of the blow-up  is made by  obtaining  the convergence of the solutions to the
$8\pi$-delta-Dirac measure concentrated at the center of mass of $n_0$.
See  also  \cite{BKLN} for a study of radially symmetric solutions for $M=8\pi$ on $\R^2$,   \cite{HV2} on a ball,
and on a domain \cite{Nagai2,Nagai3}, for instance. 
We can also refer to \cite[Th. 2 and 3, p.24-25] {Senba-Suzuki-1} for the finite number of 
 blow-up points, and for the relation between blow-up points and chemotatic collapse (Dirac measure).
\end{rem}

\begin{rem}
In the course of the proof of  Theorem \ref{tcalphazero}, we shall  see    that  the intriguing  inequality  \eqref{heatcondit2} 
is a consequence of the fact that there exists a  {\it mild}  solution $(n_t)$ of the (PKS)  system with $n_0$ as initial data. 
 More precisely,   the weak continuity at $0^+$ of $n_t$ "attaches" the free initial data $n_0$ to  the solution  $(n_t)$ for  $t>0$.
\end{rem}

\begin{rem}
The left-hand side term of the inequality  \eqref{heatcondit} is  equal to  $(4\pi T^* )p_{T^*}\star n_0(z)$.
This term is  clearly  related to the linear part of the (PKS) system  via the heat semigroup $e^{-t\Delta}$
 without  mentioning any explicit connection with the non-linear non-local term of the (PKS) system. 
The right-hand side term of \eqref{heatcondit} is  mainly related to the non-linearity of the  (PKS) system
via the rational function $R(M)$ of the mass $M$ defined in Remark \ref{rmrm}.
\end{rem}

{\bf Proof of Theorem \ref{tcalphazero}}.
The proof of  \eqref{tcnormal} is essentially a re-interpretation of the proof of the blow-up 
phenomenon when $M>8\pi$ given in \cite[p.397]{Dongyi Wei}  in  order to evaluate  blow-up time upper bounds.  
The inequality  \eqref{heatcondit2}  is due to  Dongyi Wei \cite{Dongyi Wei}.
For the sake of clarity, we  provide the main part of the arguments,
with some   comments inserted in the course of the proof. This paper is dedicated more precisely to the study of estimates of the critical time
$T_c^*(n_0)$.

 \vskip0.3cm
Fix $0\leq n_0\in L^1$, $n_0\neq 0$. 
We assume that the mass $M=\int_{\Ri^2} n_0(x)\, dx$  satisfies $M>8\pi$ and we consider  $n_t(x)=n(x,t)$ a mild solution of the (PKS) system on $\R^2$.
\vskip0.3cm

 {\bf (1)  (i)  Proof of  the inequality  \eqref{heatcondit2} and  the finiteness of $T^*$}.
\vskip0.3cm

The main tool  is the introduction of  the following  specific heat regularization function  $t\rightarrow {\tilde n}(z,t,s)=p_{s-t}\star n_t(z)$ for all  
$(t,s)$ satisfying $t<T^*$ and $t < s <+\infty$ for fixed  $z\in \R^2$. 
Here $p_t(x)=({4\pi t})^{-1}\exp(-\vert x\vert^2/4t)$ denotes the heat kernel on $\R^2$.
(At this stage,  note that $s$ is not necessarily bounded by $T^*$ itself).
The function $t\rightarrow {\tilde n}(z,t,s)$ can be seen as a  perturbation and a regularization
of the solution  $(n_t)_{t\in (0,T^*)}$ of the  (PKS) system  by the heat semigroup. 
It consists also to consider the (PKS) solution  $n_t$ at  time $t$ as the initial data (in $L^1$)  of the heat equation,
and to compute the solution of the  heat equation at  time $s-t$ for a priori  any $s>t$ by convolution with the heat kernel $p_{s-t}$.
This  construction is suggested by the non-linear term of the  Duhamel formula \eqref{duhamel}. 
 The role of reversing the time by $-t$ in the heat kernel is that the time derivative of ${\tilde n}$ 
 will  depend only on the non-linear (here bilinear) term of the  (PKS) system.  
 The linear part represented by the Laplacian disappears in this process. Indeed, we have
 $$
 \partial_t{\tilde n}(z,t,s)
 =\frac{1}{4\pi} \int_{\Ri^2} \int_{\Ri^2} K_1(x-z,y-z,s-t) n_t(x)n_t(y)\, dxdy,
 $$
 where  the  convolution kernel  in space variables is given by
 $$
 K_1(x,y,t) 
 =
 -\frac{
 \left[ \nabla p_t(x)-\nabla p_t(y)\right]. (x-y)}{\vert x-y\vert^2}.
 $$
 This kernel $K_1$ is independent of the (PKS) solution. 
 See   \cite[p.392-393]{Dongyi Wei} for the  details of computation.
 In   \cite[Prop.3.1]{Dongyi Wei}, 
the following  monotonicity inequality  for the map $t\mapsto  {\tilde n}(z,t,s)$ has been proved 
\begin{equation}\label{monoton}
\frac{3M \, {\tilde n}(z,t,s)}{8\pi (s-t)}
-\frac{M^2}{16\pi^2 (s-t)^2}
\leq
\partial_t{\tilde n}(z,t,s),
\end{equation}
for all   $z\in \R^2$ and all $(t,s)$ such that  $0<t< {\rm min}(s,T^*)$.
This result follows from the existence of a solution and the structure  of the (PKS) system,
 and some  geometric estimates on  the heat kernel  $p_t$  on $\R^2$.
Multiplying the inequality   \eqref{monoton} by the positive term $(s-t)^a$ with  $a:=\frac{3M}{8\pi}$ leads to 
$$
\partial_t\left[
(s-t)^{a}
{\tilde n}(z,t,s)\right]
=
-a(s-t)^{a-1}{\tilde n}(z,t,s)\
+ (s-t)^{a}\partial_t
{\tilde n}(z,t,s)
\geq 
\frac{-M^2}{16\pi^2}(s-t)^{a-2},
$$
for all   $z\in \R^2$ and all  $0<t<s<T^*$.
Now, by integrating  over the interval  $(0,s)$ with respect to $t$  we   obtain the following  inequality
$$
\lim_{t\rightarrow s^-} (s-t)^{a-1}(s-t) p_{s-t}\star n_t(z)
-
\lim_{t\rightarrow 0^+} (s-t)^{a}p_{s-t}\star n_t(z)
=
0 -s^a p_{s}\star n_0(z)
$$
$$
\geq 
\frac{- M^2}{16\pi^2 (a-1)}s^{a-1}
=
\frac{- M^2}{2\pi(3M-8\pi)}s^{a-1},
$$
for all $0<s<T^*$.
The first limit is zero since  we have
$$
0\leq (s-t)^a p_{s-t}\star n_t(z)\leq \frac{ (s-t)^{a-1} }{4\pi}
\int_{\Ri^2} n_t(x)\, dx=\frac{(s-t)^{a-1} }{4\pi} M,
$$ 
 for all $0<t<s<T^*$, and $a>1$.
The second limit is obtained by using  in a crucial way  the weak continuity at $t=0^+$ of  the solution $n\in C_{w}([0,T^*), L^1(\R^2))$
 which  {\em attaches} the free data  $n_0$ to the solution $(n_t)_{t\in [0, T^*)}$ of the (PKS) system.
Note also that integration on $t$  over $(0,s)$ imposes necessarily that $s<T^*$ in the conclusion since $0<t< {\rm min}(s,T^*)$.
Finally  after simplification, we  get 
$$
H_{z,n_0}(s):=4\pi s p_s\star n_0(z)
\leq
\frac{2M^2}{3M-8\pi}.
$$
Thus,  the inequality  \eqref{heatcondit2} is proved for any non-negative initial data $n_0\in L^1(\R^2)$
of the (PKS) system, and  for all $0<s<T^*$ and all $z\in \R^2$. 
\vskip0.3cm

 We now recall how we deduce  from \eqref{heatcondit2}  that  $T^*$ is finite when $M>8\pi$.
To prove this, we suppose the contrary, i.e.  $T^*=+\infty$.
 We  can take the limit as $s$ goes to infinity in \eqref{heatcondit2} and, by monotone convergence theorem, we obtain 
 $$
M= \lim_{s\rightarrow +\infty} H_{z,n_0}(s)\leq \frac{2M^2}{3M-8\pi}.
 $$
This implies $M\leq  8\pi$. 
Contradiction.
\vskip0.3cm

 {\bf(1)(ii)  Proof of  \eqref{heatcondit}}.
We have seen in (1)(i) that $T^*$ is finite.
Now by Fatou lemma, we can write
$$
 H_{z,n_0}(T^*):=\int_{\Ri^2} \exp\left(-\frac{\vert x -z\vert^2}{4T^*}\right)n_0(x)\, dx
 = \int_{\Ri^2} \liminf_{s\rightarrow  T^{*,-}}
 \exp\left(-\frac{\vert x -z\vert^2}{4s}\right)n_0(x)\, dx
 $$
 $$
 \leq
  \liminf_{s\rightarrow T^{*,-}}  \int_{\Ri^2} 
 \exp\left(-\frac{\vert x -z\vert^2}{4s}\right)n_0(x)\, dx
\leq 
\frac{2M^2}{3M-8\pi}. 
$$
This proves \eqref{heatcondit}.
\vskip0.3cm

 {\bf(1)(iii)  Proof of the converse: \eqref{heatcondit} $\Rightarrow $ \eqref{heatcondit2}}.
Conversely, the inequality  \eqref{heatcondit} implies  the inequality \eqref{heatcondit2}  
because $s\in (0, +\infty) \mapsto  H_{z,n_0}(s)$ is  clearly an increasing function.
\vskip0.3cm

{\bf (2)}
Let $z\in \R^2$, $n_0\in L^1(\R^2)$  be fixed. We also assume that $n_0\neq 0$  and  $n_0\geq 0$. It is easy to see that the function 
$$
s\in (0,+\infty) \rightarrow H_{z,n_0}(s):= \int_{\Ri^2} \exp\left(-\frac{\vert x -z\vert^2}{4s}\right)n_0(x)\, dx \in (0,M),
$$
 is continuous and strictly increasing on $(0,+\infty)$ with range $(0,M)$.
 Indeed, by monotone convergence theorem, we deduce that  $H_{z,n_0}(0^+)=0$, $H_{z,n_0}(+\infty)= M$ and  $H_{z,n_0}$ is continuous on $(0,+\infty)$.
The function $H_{z,n_0}$ is strictly increasing because $n_0\geq 0$ and $n_0\neq 0$.
Hence, $H_{z,n_0}$ is a  bijection from  $(0,+\infty)$ onto $(0,M)$. We denote by $H_{z,n_0}^{-1}:(0,M)\rightarrow (0,+\infty)$
 its inverse function. It is  also a strictly increasing continuous function.
For  $M>8\pi$,  we have  $0<\frac{2M^2}{3M-8\pi}< M$. Hence,
the value $T_{c,n_0}^*(z)=H_{z,n_0}^{-1}\left(\frac{2M^2}{3M-8\pi}\right)$  is well-defined and  finite.
By \eqref{heatcondit}, i.e.  $H_{z,n_0}(T^*)\leq \frac{2M^2}{3M-8\pi}$, and by  monotonicity of $H^{-1}_{z,n_0}$, 
this  is equivalent to the following bound 
$$
T^*\leq H_{z,n_0}^{-1}\left(\frac{2M^2}{3M-8\pi}\right)=:T_{c,z}^*(n_0) <+\infty,
$$
 for all $z\in \R^2$. Thus,  we obtain simultaneously a quantitative  upper bound on the blow-up  time $T^*$, and  again its finiteness.
Now taking the infimum over $z$, we deduce  the upper bound  \eqref{tcnormal} on the maximal existence   time $T^*$.
\vskip0.3cm

{\bf (3) (i)}
Now assume that   there exists $z_0\in\R^2$ such that we have
$H_{z,n_0}(s)\leq H_{z_0,n_0}(s)$ for all $s>0$ and all $z\in \R^2$.
 Since $H_{z,n_0}$  and $H_{z_0,n_0}$ are increasing and $H_{z,n_0}(0^+)=H_{z_0,n_0}(0^+)=0$, $H_{z,n_0}(+\infty)=H_{z_0,n_0}(+\infty)= M$,
  then we have  
 $ H^{-1}_{z_0,n_0}(m) \leq H_{z,n_0}^{-1}(m)$ for all $m\in (0,M)$ and all $z\in \R^2$. Hence, we obtain 
$$
\inf_{z\in \Ri^2} H_{z,n_0}^{-1}(m)= H^{-1}_{z_0,n_0}(m),
$$
for all $m\in (0,M)$.
Applying this result to $m=\frac{2M^2}{3M-8\pi}$ with $M>8\pi$, we obtain  
$$
T_c^*(n_0):= \inf_{z\in \Ri^2} H_{z,n_0}^{-1}\left(\frac{2M^2}{3M-8\pi}\right)
=H^{-1}_{z_0,n_0} \left(\frac{2M^2}{3M-8\pi}\right)=T_{c, z_0}^*(n_0).
$$
This concludes the first part of the assertion (3).
\vskip0.3cm

{\bf (3) (ii)}
  On $\R^n$,  let  $f$ and $g$ be two non-negative non-increasing radially symmetric 
integrable functions, then it is known that the convolution $f\star g$ is also a
non-negative non-increasing radially symmetric  integrable function. 
 In particular,  if moreover $z\mapsto  f\star g(z)$ is continuous, we deduce that $0\leq f\star g(z)\leq  f\star g(0)$,
 for all $z\in \R^n$.
Here, we have  assumed that  $g=n_0$ is a non-negative non-increasing radially symmetric  integrable  function.
Since $f=p_s$ is also a  non-negative non-increasing radially symmetric  integrable   function,
then  the function $z\in \R^2 \mapsto H_{z,n_0}(s)=4\pi s \, p_s\star n_0(z)$ is a (continuous) non-negative non-increasing radially symmetric
 integrable  function.

As a consequence, we have $H_{z,n_0}(s)\leq  H_{0,n_0}(s)$  for all $s>0$ and all $z\in \R^2$. We conclude by   applying  (3)(i) proved just above
with  $z_0=0$. This concludes the second part of the assertion (3).
\vskip0.3cm

{\bf (4)}
First note that  the masses $M:=M(n_0)=\int_{\Ri^2} n_0(x)\, dx=M(m_0)$ are equal since the Lebesgue measure is translation-invariant.
 Moreover, it is easy to check that
$H_{z,m_0}(s)=H_{z+z_0,n_0}(s)$ for all $z,z_0\in \R^2$ and all $s>0$ with $m _0(x)=n_0(x+z_0)$.
This yields
$$
T_c^*(m_0):= \inf_{z\in \Ri^2} H_{z,m_0}^{-1}(L(M))
=
\inf_{z\in \Ri^2} H_{z+z_0,n_0}^{-1}(L(M))
=\inf_{z\in \Ri^2} H_{z,n_0}^{-1}(L(M))=:
T_c^*(n_0),
$$
where  $L(M)= \frac{2M^2}{3M-8\pi}$.  
This leads to the stated  invariance.
\vskip0.3cm
This concludes the proof of Theorem \ref{tcalphazero}.
\hfill$\square$
\vskip0.5cm

In the second part of this section, we are interested in estimating the critical time 
$T^*_c=T^*_c(n_0)$ defined in Theorem \ref{tcalphazero} for  the particular case  
where $n_0$ is a  non-increasing  $z_0$-radially symmetric initial data.  
Recall that a function $n_0$ is said radially symmetric on $\R^n$  if, for all $x,y\in \R^n$ 
such that $\vert x\vert=\vert y\vert$ then $n_0(x)=n_0(y)$. We say that 
$n_0$ is non-increasing radially symmetric  on $\R^n$ if $n_0$ satisfies  for all $x,y\in \R^n$ such that 
$\vert x\vert\leq \vert y\vert$, we have  $n_0(x)\geq n_0(y)$.
This is equivalent to say  that $n_0$  is radially symmetric and the function ${\tilde n}_0$ defined by 
$r\mapsto {\tilde n}_0(r)=n_0(x)$ for $x\in \R^n$  with $r=\vert x\vert$, is  non-increasing  on $[0,+\infty)$.
Fix $z_0\in \R^n$, we say that a function $n_0$ is a $z_0$-radially symmetric  (resp. non-increasing $z_0$-radially symmetric)  
function if  $x\mapsto m_0(x):=n_0(x+z_0)$ is a  radially symmetric (resp. non-increasing radially symmetric) function is the sense given above.
\vskip0.3cm

In the following results, we shall use the Laplace  transform of a function $f$ denoted  by
${\mathcal L}f (v)=\int_0^{+\infty} e^{-v u}f(u)\, du$ ($v> 0$). 
This  Laplace  transform  will be used in association  with  radially symmetric initial data $n_0$.

\begin{theo}\label{tclaplace}
Under the assumptions of Theorem \ref{tcalphazero}
and denoting by $T^*_c(n_0)$ the critical time defined by \eqref{tcnormal}. We have the next statements.
\begin{enumerate}
\item
If $n_0$  is a $z_0$-radially symmetric non-negative integrable  function then
\begin{equation}\label{laplacegal2}
T^*\leq  T^*_c(n_0)\leq 
\frac{1}{ 4({\mathcal L} f)
^{-1}\left(\frac{L(M)}{\pi}\right)
},
\end{equation}
with   $f(u)={\tilde m_0}(\sqrt{u})$, $u>0$, where ${\tilde m_0}(\vert x\vert)=m_0(x)=n_0(x+z_0)$, $x,z_0\in \R^2$, and $L(M)=\frac{2M^2}{3M-8\pi}$.
\item
Moreover, suppose that
$n_0$  is a  non-increasing $z_0$-radially symmetric non-negative integrable   function. 
Then we have
\begin{equation}\label{laplacegal1}
T^*\leq T^*_c(n_0)=  
\frac{1}{ 4({\mathcal L} f)
^{-1}\left(\frac{L(M)}{\pi}\right)
},
\end{equation}
with $f$ as in \eqref{laplacegal2}.
\end{enumerate}
\end{theo}

Note that the upper bound on  $T^*_c(n_0)$ in  \eqref{laplacegal2} is now an equality
in  \eqref{laplacegal1} when $n_0$  is a  non-increasing $z_0$-radially symmetric  function.
So, the estimate   \eqref{laplacegal2}  of $T^*_c(n_0)$ is sharp for the  sub-class of non-increasing $z_0$-radially symmetric data $n_0$.
\vskip0.3cm 

{\bf Proof}. 
By applying  \eqref{tcnormal} and the statement (4)   of Theorem \ref{tcalphazero}, we obtain 
$$
T^*\leq 
T_c^*(n_0)=T_c^*(m_0)
\leq  T_{c,0}^*(m_0)
=H_{0,m_0}^{-1}\left(L(M)\right).
$$
Now, since the function $m_0$ is a radially symmetric, the function  $H_{0,m_0}$ can expressed via  the Laplace transform as follows,
$$
H_{0,m_0}(s)=
\int_{\Ri^2} \exp\left(-\frac{\vert x  \vert^2}{4s}\right)m_0(x)\, dx
=
\int_{\Ri^2} \exp\left(-\frac{\vert x  \vert^2}{4s}\right){\tilde m_0}(\vert x\vert )\, dx.
$$ 
By the change of variables  using polar coordinates, we deduce that
$$
H_{0,m_0}(s)=
2\pi \int_0^{+\infty} e^{-r^2/4s} {\tilde m_0}(r)r dr
=
\pi \int_0^{+\infty} e^{-r^2/4s}  {\tilde m_0}(\sqrt{r^2})\, d(r^2)
$$
$$
=
\pi \int_0^{+\infty} e^{-\rho/4s}   {\tilde m_0}(\sqrt{\rho})\, d\rho
=
\pi \int_0^{+\infty} e^{-\rho/4s}f(\rho)\, d\rho
=
\pi ({\mathcal L} f)(1/4s), 
$$
with $f(\rho)= {\tilde m_0}(\sqrt{\rho})$. The last line is obtained by the change of variable $\rho=r^2$.
The  function  $v\rightarrow {\mathcal L}f(v)$ is invertible since $ H_{0,m_0}$ is invertible.
Finally, we conclude that
$$
H^{-1}_{0,m_0}(L(M))=\frac{1}{ 4({\mathcal L}f )
^{-1}\left(\frac{L(M)}{\pi}\right)
}
.
$$
This proves the inequality \eqref{laplacegal2}. 
\vskip0.3cm
We now prove  the inequality \eqref{laplacegal1} as follows.
Since  $n_0$  is a  non-increasing $z_0$-radially symmetric  function, 
then $m_0$ is a non-increasing radially symmetric  function.
We  apply (3)  of Theorem \ref{tcalphazero}  to  obtain the next equality
$$ 
T_c^*(n_0)=T_c^*(m_0)=\inf_{z\in \Ri^2} H_{z,m_0}^{-1}(L(M)) 
=
H^{-1}_{0,m_0}(L(M))=
\frac{1}{ 4({\mathcal L}f )
^{-1}\left(\frac{L(M)}{\pi}\right)
}.
$$
\vskip0.3cm
This concludes the proof of Theorem \ref{tclaplace}.
\hfill $\square$
\vskip0.3cm

\noindent
The other sections of this paper  are organized  as follows.
\vskip0.3cm

\indent
In  Section \ref{tcsection}, since  generally  it is  difficult to  compute  $T_c^*(n_0)$ explicitly,  we extend our investigation of upper bounds
on $ T^*$  by  providing  various (more or less explicit) upper bounds  on the bound $T_c^*(n_0)$ itself.
We also prove lower bounds on $T_c^*(n_0)$ in Section \ref{lowerboundsect} showing (in some cases)  
the sharpness of the corresponding upper bounds on $T_c^*(n_0)$. 
In Section \ref{examples}, we finally exhibit  explicit  bounds on $T_c^*(n_0)$  for several  families of examples of initial data.
In particular,  we can  compute explicitly $T_c^*(n_0)$ for some examples of radially symmetric initial data $n_0$
by applying  Theorem \ref{tclaplace} using the Laplace transform.

\section{Estimates on the critical time bound   $T^*_c(n_0)$}\label{tcsection}
In  general, the functions $H_{z,n_0}^{-1}$ of  Theorem \ref{tcalphazero} are not given  explicitly.
Nevertheless, from this theorem  we are able  to derive various types of  explicit upper bounds on $T^*$ 
from  those on  $T^*_c(n_0)$,  for all $0\leq n_0\in L^1(\R^2)$, with or without any specific additional assumptions.
Such  estimates may be useful by  providing an explicit interval  of time $(0,\tilde{T})\supseteq  (0,T^*_c(n_0)) \supseteq  (0,T^*)$ 
 for numerical simulations on which the blow-up will certainly be observed. 

We  first obtain an upper bound on $T^*(n_0)$ as an infimum over two parameters of a function 
expressed in terms of  convolutions with Gaussian-type weights of the initial data  $n_0$.
We denote by $\ln_+ v=\sup(\ln v, 0)$,  $v>0$, and  
$\vert\vert f \vert\vert_{\infty}=\sup_{ z\in \Ri^2}\vert f(z)\vert $ for a bounded function $f$.

\begin{cor}\label{corapplica0}
Let $T^*=T^*(n_0)$ be the maximal existence  time of a solution of the (PKS) system  with initial data $0\leq n_0\in L^1$
satisfying $M=\vert\vert n_0\vert\vert_{1}>8\pi$. Let $L(M)=\frac{2M^2}{3M-8\pi}$.
Then we have
\begin{equation}\label{app0}
  T^* \leq T^*_c(n_0)
   \leq
T^{*}_{c_1}:=
\inf_{q>1,\, \lambda>0} 
\lambda q^{-1/q}
\left[\ln_+\left(  \frac{
\vert\vert  \omega_{q,\lambda}   \star n_0 \vert\vert_{\infty}} {L(M)} \right)
\right]^{-1/q},
\end{equation}
where $\omega_{q,\lambda}(x)=\exp(-c_{q,\lambda}
\vert x\vert^{\frac{2q}{q-1}}) $
with  $c_{q,\lambda}=\frac{q-1}{q} (4\lambda)^{-\frac{q}{q-1}}$,  for all $1<q<\infty$ and  all $\lambda >0$.
In particular, if $n_0$ is a non-increasing  radially symmetric function, then
$$
\vert\vert  \omega_{q,\lambda}   \star n_0 \vert\vert_{\infty}
=  \omega_{q,\lambda}   \star n_0 (0)=
 \int_{\Ri^2} \exp\left(-c_{q,\lambda}
\vert x \vert^{\frac{2q}{q-1}}\right) n_0(x)\, dx.
$$
\end{cor}

We make some comments on these results before giving the proof.
\begin{rem}
 Note that the weight $\omega_{q,\lambda}$ does not depend on $T^*$.
 In fact, the parameter $\lambda$ in the  formula \eqref{app0} plays a similar role to $T^*$ 
 but now it can be taken freely  in the full interval $(0,+\infty)$.
\end{rem}

\begin{rem}
 The value $T^{*}_{c_1}$ is  finite. We prove it as follows.
Let $p=\frac{q}{q-1}$  with  $1<q<\infty$.
For all $1<q<+\infty$ and all $z\in \R^2$, we have
by dominated convergence theorem, 
$$
\lim_{\lambda\rightarrow +\infty}
\frac{\omega_{q,\lambda}\star n_0(z)}{L(M)}
=
\frac{1}{L(M)}
\int_{\Ri^2} \lim_{\lambda\rightarrow +\infty}
 \exp\left[-\frac{1}{p} 
\left(\frac{\vert x-z\vert^{2}}{4\lambda}\right)^p\, \right]
n_0(x)\, dx
=\frac{M}{L(M)}>1,
$$
whenever $M>8\pi$.
Thus, we get for all  $1<q<+\infty$ and $\lambda \geq \lambda_q$ (large enough),
$$
0 =\ln_+1<
\ln_+\
\left( 
\frac{
\vert\vert \omega_{q,\lambda}\star n_0 \vert\vert_{\infty}
}
{L(M)}
\right).
$$
This  implies  that $T^{*}_{c_1}$ defined by \eqref{app0}  is finite.
\end{rem}

{\bf Proof}.
The main  idea  of proof is to obtain a lower bound for  the integral  $H_{z,n_0}(s)$ by a product of two functions, one  depending on  the time variable $s$
and the other one depending on  the space variable $z$, in order  to get  a lower bound of the following  form,
$$
g(s)(\omega\star n_0)(z)\leq 
 \int_{\Ri^2} \exp\left(-\frac{\vert x -z\vert^2}{4s}\right)n_0(x)\, dx
 =
 H_{z,n_0}(s), \, z\in \R^2, \quad s>0.
$$
Here, we expect that $s\mapsto g(s)$ will be  an   invertible function explicitly given,
 and the function   $x \in \R^2 \mapsto\omega(x)$ will   also  be an explicit  (bounded) weight of the space variable $x$.
\vskip0.3cm

More explicitly, we prove the corollary a  follows.
We apply the  next well-known  Young's  inequality,
$$
ab\leq \frac{a^p}{p\lambda^p}+
 \frac{\lambda^q b^q}{q},
 $$
for all  $a,b\geq 0$, $\lambda >0$,  where  $1<p,q<+\infty$  satisfies  $1/p+1/q=1$ (i.e., $p=\frac{q}{q-1}$). 
We  set   $a={\vert x-z\vert^2}/{4}$ and  $b=1/s$ 
for all $x,z\in \R^2$ and all $s>0$. Then  we deduce that
$$
\frac{\vert x-z\vert^2}{4s}
\leq 
\ell  \vert x-z\vert^{2p}
+
 \frac{k}{s^q},
 $$
 with $\ell=p^{-1}(4\lambda)^{-p}$ and $k=q^{-1}\lambda^q$.
This implies that
 $$
e^{- ks^{-q}}
 \int_{\Ri^2} 
 \exp\left(- \ell \vert x-z\vert^{2p} \right)n_0(x)\, dx
 \leq
 \int_{\Ri^2} 
 \exp\left(-\frac{\vert x-z\vert^{2}}{4s}\right) n_0(x)\, dx,
 $$
 for all $z\in \R^2$ and all $s>0$.
 By setting $\omega_{q,\lambda}(x)=\exp\left(-\ell \vert x\vert^{2p} \right)$, our result reads as a convolution inequality,
$$
e^{- ks^{-q}}
\omega_{q,\lambda}  \star n_0(z)
\leq H_{z,n_0}(s).
 $$
 Recall that $T^*_c=T^*_c(n_0):=\inf_{z\in \Ri^2} T^*_{c,z}(n_0)$, and that $T^*_{c,z}:=T^*_{c,z}(n_0)$ is defined as the unique solution of  $H_{z,n_0}(T^*_{c,z})
=L(M)$ for any fixed $z\in \R^2$.
 By monotonicity of the functions  $s\mapsto  H_{z,n_0}(s)$
 and using  the definition of $T^*_{c}$, we deduce from this convolution inequality  that
 $$
e^{- ks^{-q}}
\omega_{q,\lambda} \star n_0(z)
\leq  H_{z,n_0}(T^*_c)
\leq H_{z,n_0}(T^*_{c,z})
=L(M).
 $$
 for all  $0<s\leq T^*_c$ and all $z\in \R^2$. 
Now taking the supremum over $z\in \R^2$ and $s=T^*_c$, we obtain the next  inequality,
$$
e^{- k  (T^*_c)^{-q}}
\vert\vert \omega_{q,\lambda}  \star n_0 \vert\vert_{\infty}
\leq L(M).
$$
From the inequation \eqref{tcnormal} of Theorem \ref{tcalphazero},  
we deduce that
$$
  T^* \leq T^*_c(n_0)
\leq
k^{1/q}
\left[\ln_+\left(  \frac{
\vert\vert  \omega_{q,\lambda}   \star n_0 \vert\vert_{\infty}} {L(M)} \right)
\right]^{-1/q}
$$
$$
=
\lambda q^{-1/q}
\left[\ln_+\left(  \frac{
\vert\vert  \omega_{q,\lambda}   \star n_0 \vert\vert_{\infty}} {L(M)} \right)
\right]^{-1/q},
$$
for all $1<q<+\infty$ and  all $\lambda>0$.

Finally, by taking the infimum over the parameters  $q>1$ and $\lambda>0$, we obtain the inequality  \eqref{app0}.
\vskip0.3cm

Now we prove the last statement of Corollary \ref{corapplica0} as follows.
By assumption, $n_0$ is a  non-increasing radially symmetric  function, and since the function
$$
x\mapsto \omega_{q,\lambda} (x)=
\exp(-c_{q,\lambda}
\vert x\vert^{\frac{2q}{q-1}})
$$
has clearly  the same property,
then the convolution of these two functions, i.e. 
 $z\mapsto \omega_{q,\lambda}  \star n_0(z)$,
 is also a   non-increasing radially symmetric (continuous) function. Thus, it attains its supremum at $z=0$.
This concludes the proof of Corollary \ref{corapplica0}.
\hfill$\square$
\vskip0.3cm

Unfortunately, despite its theoretical interest, it is difficult to estimate the infimum $T_{c_1}^*$ in \eqref{app0} of 
Corollary \ref{corapplica0}, even for the simplest case of  the characteristic function of a disk, i.e.  $n_0=1_{B(0,R)}$. 
Nevertheless, in the general situation we can choose any fixed couple of parameter $q>1$ and $\lambda>0$ to bound $T_c^*(n_0)$, hence $T^*$ too.
In the next corollary, we  propose several  approaches to obtain more explicit bounds on $T_c^*(n_0)$. 
The first one is given  in terms of averages  of $n_0$ on the family of  disks  $B(z,\rho)$ of radius $\rho>0$  centered at  $z\in \R^2$, 
and valid for any  initial data $n_0$.
Next, in this corollary, we rewrite this result differently. 
The last result  mentioned in this corollary deals with the case  of data with compact  support.
All these  results are  applied to obtain explicit upper bounds on $T^*$ for several families of examples, see Section \ref{examples}.

\begin{cor}\label{corapplica1}
 Let  $T^*$ be the maximal  existence  time  of  a solution of the  (PKS) system  \eqref{pks1}-\eqref{init} with initial data $n_0\geq 0 $ of mass
$M>8\pi$. Let $L(M)=\frac{2M^2}{3M-8\pi}$.
\begin{enumerate}
\item
We have $T^*\leq T_{c}^*:=T_{c}^*(n_0)$ with
\begin{equation}\label{tdeuxc}
T_c^*\leq T_{c_2}^*:=
\inf_{\rho>0,z\in \Ri^2}
\frac{\rho^2}{4}\left[\ln_+\left(\frac{{\mathcal M}_{z}(\rho)}{L(M)}\right)\right]^{-1},
\end{equation}
where 
${\mathcal M}_{z}(\rho):=  
\int_{B(z,\rho)} n_0(x)\, dx$ and 
$B(z,\rho)$ is the Euclidean disk of radius $\rho>0$  centered at  $z\in \R^2$.
Moreover, if $n_0$ is a non-increasing radially symmetric (integrable)  function, then we have
\begin{equation}\label{radial}
T^* \leq T_c^*\leq T_{c_2}^*=
\inf_{\rho>0}
\frac{\rho^2}{4}\left[\ln_+\left(\frac{{\mathcal M}_{0}(\rho)}{L(M)}\right)\right]^{-1}.
\end{equation}
\vskip0.3cm

\item
For  each fixed  $z\in \R^2$ and $\rho>0$, let $g_z(\rho)=\frac{1}{M} \int_{B(z,\rho)} n_0(x)\, dx$. We denote by $g_z^{\leftarrow}$ 
the generalized inverse of $g$ defined by
$$
g_z^{\leftarrow}(m):=\inf\{ \rho >0 : g_z(\rho)\geq m\},\quad m\in (0,1).
$$
For all $z\in \R^2$, let 
\begin{equation}\label{atheta}
T^*_{c_2}(z):=\inf_{\rho>0}
\frac{\rho^2}{4}\left[\ln_+\left(\frac{{\mathcal M}_{z}(\rho)}{L(M)}\right)\right]^{-1}.
$$
Then we have
$$
T^*_{c_2}(z)=
\frac{1}
{
4\ln\left(
1+\frac{M-8\pi}{2M}
\right)
}.
\inf_{\theta \in (0,1)}
\left(
\frac{ 
\left[
g_z^{\leftarrow}(a^{\theta})
\right]^2
}{(1-\theta)}
\right),
\end{equation}
with $a:=\frac{L(M)}{M}=\frac{2M}{3M-8\pi} \in (0,1)$.
Hence, we also have
\begin{equation}\label{tdeuxd}
T^* \leq  T_c^* \leq T_{c_2}^*=
\frac{1}
{
4\ln\left(
1+\frac{M-8\pi}{2M}
\right)
}.
\inf_{z\in \Ri^2,\theta \in (0,1)}
\left(
\frac{ 
\left[
g_z^{\leftarrow}(a^{\theta})
\right]^2
}{(1-\theta)}
\right).
\end{equation}\vskip0.3cm
\item
Assume that $n_0$ has  compact support denoted by $K$. 
Let $i_K(z)=\sup_{x\in K}\vert x-z\vert$ for all $z\in \R^2$,
and $R_0=\inf_{z\in\Ri^2} i_K(z)$. 
Then we have
\begin{equation}\label{tctrois}
T^*\leq  T^*_{c}  \leq T_{c_3}^*:=
\frac{R_0^2}{4\ln\left(
1+\frac{M-8\pi}{2M}
\right)}.
\end{equation}
 The value $R_0$ is  also the radius of the smallest closed disk containing $K$.
 In particular, we have
\begin{equation}\label{tctroisbis}
T^* \leq  T^*_{c}  \leq  T_{c_3}^*\leq 
\frac{D^2}{12\ln\left(
1+\frac{M-8\pi}{2M}
\right)}, 
\end{equation}
where $D$ is the diameter of $K$.
\end{enumerate}
\end{cor}

Before giving the proof of this corollary, we make some comments on these results.

\begin{rem}
 By translation  invariance of $T_c^*(n_0)$ as mentioned in Theorem \ref{tcalphazero} (4),
the inequality  \eqref{radial} is still valid  when $n_0$ is a  non-increasing  $z_0$-radially symmetric function  changing 
${\mathcal M}_0(\rho)$ by ${\mathcal M}_{z_0}(\rho)$.
\end{rem}

\begin{rem}
 The ratio $\frac{1}{4\ln\left(1+\frac{M-8\pi}{2M}\right)}$
appears naturally in the proof of  \eqref{tctrois}. 
The inequality  \eqref{tdeuxd} proposes a more  general version where  this ratio appears
without  any compact support condition on the initial data $n_0$. 
\end{rem}

\begin{rem}\label{remgz}
\begin{enumerate}
\item
 The range of the non-decreasing continuous function $[0,+\infty) \ni \rho \mapsto g_z(\rho)$ is  either  the interval $[0,1)$ or $[0,1]$.
Indeed, we have $g_z(0)=0$, $\lim_{\rho\rightarrow +\infty} g_z(\rho)=1$ and $g_z$ is continuous on  its domain $[0,+\infty)$  for all $z\in \R^2$.
This function may not be a strictly increasing function on $[0,+\infty)$. In particular, it can be  constant on an interval.
For instance  if $g_z(\rho_1)=g_z(\rho_2)$ for some $0<\rho_1<\rho_2\leq +\infty$ and for some  fixed $z\in \R^2$,
then  the  (non-negative) initial data  $n_0$  is almost everywhere zero in the annulus 
$B(z,\rho_2)\setminus B(z,\rho_1)$.  It is easy to construct such examples of initial data. 
For instance, we can simply  consider the characteristic function of the disk $B(z_0,\rho_1)$ and the function $g_{z_0}$ for a fixed $z_0$
 (here, $\rho_2=+\infty$). More generally, we can consider an initial data with compact support.
\item
If $g_z$  is  strictly increasing on $[0,+\infty)$, then its range is $[0,1)$. Thus, the function
$g_z:[0,+\infty)\rightarrow [0,1)$ is a bijection and 
$g_z^{\leftarrow}$ is the usual inverse function  $g_z^{-1}:[0,1)\rightarrow [0,+\infty)$ of $g_z$.  
Indeed, the value $1$ cannot be  in the range of $g_z$, if not,
there exists $0<\rho_0<+\infty$ such that 
$g_z(\rho_0)=1$,  then we have $g_z(\rho)=1$ for all $\rho\geq \rho_0$, since $g_z$ is non-decreasing and bounded by $1$.
 Contradiction:  $g_z$ is not strictly increasing.
\item
 If the value $1$ is  in the range of $g_z$, 
 i.e.  there exists $0<\rho_0<+\infty$ such that $g_z(\rho_0)=1$,  then we deduce that $g_z^{\leftarrow}(1)$ is defined, and
$$
g_z^{\leftarrow}(1) 
:=\inf \{\rho>0, g_z(\rho)\geq 1\}
=
\inf \{\rho>0, g_z(\rho)=1\}\leq \rho_0.
$$
So, $g_z^{\leftarrow}(1)$ is defined and finite.
In fact,  $g_z^{\leftarrow}(1)$ is defined and finite if and only if  the (essential)  support of $n_0\geq 0$ is included in a closed disk 
$\overline{B(z,\rho_0})$ for some $0<\rho_0<+\infty$.
In that situation, we shall  simply say that  the initial data $n_0$ has compact support.
It remains to prove the converse. Assume that $g_z^{\leftarrow}(1)$ is defined and finite.
Let $\rho_1:=g_z^{\leftarrow}(1)<+\infty$.
Then, we have $g_z(\rho_0)=1$ for at least one  $ \rho_0\geq  \rho_1$. Thus, we obtain
 $$
0=1-g_z(\rho_0)=\frac{1}{M}\int_{B^c(z,\rho_0)} n_0(x)\, dx.
 $$
Since $n_0$ is non-negative, we deduce that $n_0=0$ almost everywhere on $\R^2\setminus  B(z,\rho_0)$, which proves the assertion.
 \item
If  the value $1$ is not in the range of $g_z$ (i.e. $n_0$ has no compact support), and if needed, we  shall let  $g_z^{\leftarrow}(1)=+\infty$.
\item
The function $g_z^{\leftarrow}$ will be  called the  {\em generalized  right-inverse} of $g_z$.
It   can also be  called the radial cumulative distribution function centered at $z\in \R^2$
of the density of probability
$d\mu(x)=\frac{1}{M} n_0(x)\, dx$.
\end{enumerate}
\end{rem}

\begin{rem}
During  the proof, we shall see  that both inequalities  \eqref{radial} and \eqref{tdeuxd} are equivalent. 
 But  in  the  second  inequality   \eqref{tdeuxd},  the following   quantity 
 $$
 \frac{1} {4\ln\left(1+\frac{M-8\pi}{2M}\right)}
 $$
appears naturally  in the  estimate of the upper bound of $T^*_c(n_0)$, and also
 in \eqref{tctrois} and  \eqref{tctroisbis}.
Note that this term is independent of the particular  shape of the initial data $n_0$, and depends only on the mass $M$.
 Such a logarithmic  term   is mentioned here because it also occurs in Theorem \ref{lowe1} below for  lower bounds to $T^*_c(n_0)$.
\end{rem}

\begin{rem}\label{jungth}
The quantity $i_K(z)=\sup_{x\in K} \vert x-z\vert=\vert x_0-z\vert$ is the distance from $z$ to one of the farthest  point  $x_0=x_0(z,K)$ in $K$. 
Such a point  $x_0$  exists  because the function $x\rightarrow \vert x-z\vert$ is continuous on the compact  set  $K$. 
Hence, the supremum is attained and is, in fact, a maximum. In particular, we have for all $z\in \R^2$,
 $$
 K\subset {B^{\,\prime}(z,i_K(z))}
 : =
 \{y\in \R^2: \vert z-y\vert\leq i_K(z) \},
 $$ and $i_K(z)$ is exactly the smallest radius $R\geq 0$ such that 
 $K\subset {B^{\,\prime}(z,R)}$ for this fixed $z\in \R^2$.
Now, the quantity $R_0=\inf_{z\in\Ri^2} i_K(z)$ consists  in choosing ideally  a closed disk  of minimal radius, i.e. $R_0$,  containing $K$. 
More precisely, it can be proved that $z\rightarrow i_K(z)$ is $1$-Lipschitz (so, it is continuous), and 
$$
\lim_{\vert z\vert\rightarrow+\infty}i_K(z)=+\infty.
$$
Since $i_K\geq 0$ then the infimum defining $R_0$ is attained at some point $z_0\in \R^2$, i.e.
$R_0=i_K(z_0)=\vert z_0-x_0\vert$ for some $x_0\in K$. Then 
$K\subset  {B^{\,\prime}(z_0,R_0)}$ where $R_0$ is easily seen as  the minimal radius $R$  such that
$K$ is enclosed in a closed disk  ${B^{\,\prime}(z,R)}$ for some $z\in \R^2$ ($z$ is not necessarily in $K$).
See the examples of initial data $n_0$ given by the characteristic function of $K$ where $K$ is a closed disk or an annulus 
treated in Section \ref{examples}.
 In any metric space, it can be  easily shown that 
 $$
 \frac{1}{2} D \leq R_0\leq D,
 $$
 where $D:={\rm diam} (K)=\sup_{x,y\in K} \vert x-y\vert$ is the diameter of the  bounded (compact) set $K$.
 In $\R^2$, it may happen that $\frac{1}{2}  D <R_0$ ($R_0=D/\sqrt{3}$ for the equilateral triangle) or $R_0 < D$ (always true in $\R^n$).
 Indeed, by Jung's Theorem \cite{wiki} we  have on $\R^n$,
 $$
 R_0\leq D\sqrt{\frac{n}{2(n+1)}}.
 $$
 The case of equality is attained by the regular $n$-simplex. 
In particular on  $\R^2$, we obtain
 $$
R_0\leq D/\sqrt{3}.
$$
\end{rem}
 
\begin{rem}
 The upper bound $T^*_c(n_0)$ of $T^*$ depends not only of the mass of the initial  data  $n_0$ but also of its {\em shape}.
Indeed,  for instance,  the radial distribution ${\mathcal M}_z(\rho),\rho>0, z\in \R^2$  appearing in (1) of Corollary \ref{corapplica1} "encodes"
a part of the information on the {\em shape} of $n_0$.
In Corollary  \ref{corapplica1} (2), the geometric information is encoded by the generalized inverse 
$g^{\leftarrow}_z, z\in \R ^2$. 
For explicit {\em  shape} expressions, see the examples  described  in Section \ref{examples}.
\end{rem}

\begin{rem}
 It can be noticed that the result (1) of Corollary  \ref{corapplica1} is sharp in the sense that the set of  information 
$\{{\mathcal M}_z(\rho),\rho>0, z\in \Ri^2\}$ is equivalent to the knowledge of the initial data $n_0$ since $n_0\in L_{loc}^1$. 
Indeed, by  Lebesgue differentiation theorem, we have
$$
\lim_{\rho\rightarrow 0^+} \frac{1}{\vert B(z,\rho)\vert} {\mathcal M}_z(\rho)=n_0(z),
$$
for almost all  $z\in \R^2$.  
\end{rem}
  
{\bf Proof of Corollary \ref{corapplica1}}
\vskip0.3cm
(1) (a)
Our proof starts from \eqref{heatcondit} of Theorem \ref{tcalphazero}. 
To simplify the notation, we set $T=T_c^*(n_0)$. 
Since $n_0\geq 0$, we have 
for all $z\in \R^2$ and all $\rho>0$,
$$
e^{-\frac{\rho^2}{4T}} {\mathcal M}_{z}(\rho)
=
\int_{B(z,\rho)}  e^{-\frac{\rho^2}{4T} } n_0(x)\, dx
\leq
 \int_{\Ri^2} \exp\left(-\frac{\vert x -z\vert^2}{4T}\right)n_0(x)\, dx
=
L(M),
$$
with $L(M)=\frac{2M^2}{3M-8\pi}$.
It follows easily that, for all $z\in \R^2$ and  all $\rho>0$, 
$$
T\leq
\frac{\rho^2}
{
4\ln_+\left( 
\frac{{\mathcal M}_{z}(\rho)}{L(M)}\right)
}
\leq +\infty
.
$$
Now by taking the infimum over $z\in \R^2$ and  $\rho>0$, we deduce that
$$
T=T_c^*(n_0) \leq T_{c_2}^*:=
\inf_{\rho>0, z\in \Ri^2}
\frac{\rho^2}{4}\left[\ln_+\left(\frac{{\mathcal M}_{z}(\rho)}{L(M)}\right)\right]^{-1}.
$$
It remains to prove that  $T^*_{c_2}$ is finite for any non-negative integrable function $n_0$. 
The proof is as follows. We note that  for all $z\in \R^2$,
$$ 
\lim_{\rho\rightarrow+\infty} 
\ln_+\left(\frac{{\mathcal M}_{z}(\rho)}{L(M)}\right)
=
\ln_+\left(\frac{M}{L(M)}\right)
>0,
$$ by continuity of $\ln_+$ and monotone convergence theorem.
The last inequality is due to the fact that $L(M)<M$ when $M>8\pi$.
Then for one (all)  $z\in \R^2$, there exists $0<\rho_0<+\infty$ large enough such that,
$$
\ln_+\left(\frac{{\mathcal M}_{z}(\rho_0)}{L(M)}\right)>0.
$$
Hence, it yields
 $$
0\leq 
\frac{\rho_0^2}{4}\left[\ln_+\left(\frac{{\mathcal M}_{z}(\rho_0)}{L(M)}\right)\right]^{-1}< +\infty.
$$
This implies the finiteness of $T^*_{c_2}$.
The upper bound  \eqref{tdeuxc}  on $T_c^*(n_0)$ is now proved. 
\vskip0.3cm

(1) (b) For the proof of the second part  of (1), we just need to prove for any  non-negative non-increasing radially symmetric function $n_0$ that
\begin{equation}\label{sup}
\sup_{z\in \Ri^2} 
{\mathcal M}_{z}(\rho)
={\mathcal M}_{0}(\rho), \quad \rho>0.
\end{equation}
This result is well-known but we provide below  details for completeness.
Let $z\in \R^2$ be fixed.
For all $x\in  B(0,\rho)\setminus U$ and all $y\in  B(z,\rho)\setminus U$ with $U=B(0,\rho)\cap B(z,\rho)$,
we have  $\vert x\vert \leq \rho\leq \vert y\vert$. This implies that $n_0(x)\geq n_0(y)$
because   $n_0$ is a non-increasing radially symmetric function.
We integrate this last inequality with respect to  the couple of variables $(x,y)\in (B(0,\rho)\setminus U)\times (B(z,\rho)\setminus U)$, and get by Fubini's theorem,
$$
\vert 
B(z,\rho)\setminus U
\vert
.
\int_{
B(0,\rho)\setminus U} n_0(x)\,dx
\geq
\vert 
B(0,\rho)\setminus U
\vert .
\int_{
B(z,\rho)\setminus U} 
n_0(y)\,dy.
$$
Here,  $\vert A\vert$ denotes  the Lebesgue measure of the measurable set $A\subset \R^2$.
Since $U \subset B(z,\rho)$ and  $U \subset B(0,\rho)$, we have
$$
\vert  B(z,\rho)\setminus U \vert
=\vert  B(z,\rho)\vert -\vert U \vert
=\vert  B(0,\rho)\vert -\vert U \vert
=
\vert  B(0,\rho)\setminus U\vert
.
$$
It immediately yields 

\begin{equation}\label{boule}
\int_{
B(0,\rho)\setminus U} n_0(x)\,dx
\geq
\int_{
B(z,\rho)\setminus U} 
n_0(y)\,dy
\end{equation}

when $\vert B(z,\rho)\setminus U \vert>0$.
If $\vert B(z,\rho)\setminus U \vert=0$ then  $B(z,\rho)= B(0,\rho)\cap B(z,\rho)\subset B(0,\rho)$. This implies that $z=0\in \R^2$
and   \eqref{boule} holds trivially.
Thus,  we deduce that
$$
{\mathcal M}_{z}(\rho)
-
{\mathcal M}_{0}(\rho)
=
\int_{
B(z,\rho)\setminus U} 
n_0(y)\,dy
-
\int_{
B(0,\rho)\setminus U} 
n_0(x)\,dx
\leq 0.
$$
Finally, we have for all $z\in \R^2$ and all $\rho>0$,
$$
{\mathcal M}_z(\rho)\leq {\mathcal M}_0(\rho)
.
$$
Thus,  the formula  \eqref{sup} follows.
\vskip0.3cm

(2) The proof relies on a change of variables of exponential type.
We  first start with   some useful properties of the distribution function $g_z$ for any $z\in \R^2$,
see also Remark \ref{remgz} above.
We denote $g_z$  by $g$ at some places  for short.

The function $g:[0,+\infty)\rightarrow [0,1]$ is a non-decreasing continuous function with range included in $[0,1]$.
Since $g(0)=0$, $g(+\infty)= 1$ and  $g$ is continuous, then  for each $m\in (0,1)$ 
 there exists  $\rho_1>0$ such that $g(\rho_1)=m$ by the intermediate value theorem.
Let $G(m)=\{ \rho>0: g(\rho)=m\}$ with $m\in (0,1)$.
Then the set  $G_m$  is not empty and bounded from below by $0$, hence the infimum of  $G_m$ exists for all $m\in (0,1)$.
It is  also a minimum because  $G(m)$ is closed ($g$ is continuous).  
We set  $g^{\leftarrow}(m):={\rm min\,} G_m$. The function $g^{\leftarrow}$ is  a right-inverse of the function $g$ since we have
$g\circ g^{\leftarrow}(m)=m$ for all $m\in (0,1)$. This is simply  due to the fact that $g^{\leftarrow}(m)\in G(m)$.
Note that $g^{\leftarrow}$ is not necessarily a left-inverse of the  function $g$, in particular if $g_z$ is constant on some (closed) interval. 
Moreover, it can be shown by  monotonicity of $g$ that we also have
$$
g^{\leftarrow}(m)=\inf G_m =\inf\{ \rho>0: g(\rho)\geq m\}.
$$
Hence, the function $g^{\leftarrow}$ is non-decreasing on $(0,1)$ but not necessarily continuous.
Its range is included in $(0,+\infty)$.
Of course, if $g$ is strictly increasing  from  $(0,+\infty)$ onto $(0,1)$ then $g^{\leftarrow}=g^{-1}$ is the usual inverse function to $g$.
\vskip0.3cm

Now we are in position to prove  \eqref{tctrois} of Corollary \ref{corapplica1}.
Let $z\in \R^2$ be fixed. We set $a=\frac{L(M)}{M}=\frac{2M}{3M-8\pi}$.
Since $M>8\pi$, we have $a\in (0,1)$. By definition of  $T^*_{c_2}(z)$, we can write for all $\rho>0$,
$$
T_{c_2}^*(z)
\leq 
\frac{\rho^2}{4}\left[\ln_+\left(\frac{{\mathcal M}_{z}(\rho)}{L(M)}\right)\right]^{-1}
=
\frac{\rho^2}{4}\left[\ln_+\left(\frac{g_z(\rho)}{a}\right)\right]^{-1},
$$
because 
$ {\mathcal M}_{z}(\rho)/L(M)=g_z(\rho)/a$.
Let $\theta\in (0,1)$ then we have $a^{\theta}\in (0,1)$. 
We set $\rho_*=g_{z}^{\leftarrow}(a^{\theta})$ then we have $\rho_*\in (0,+\infty)$ and $g_z(\rho_*)=a^{\theta}$.
The inequality just above with $\rho=\rho_*$ implies that
$$
T^*_{c_2}\leq 
\frac{\left[ g_z^{\leftarrow}(a^{\theta})\right]^2}
{
4\ln_+\left(a^{\theta-1}\right)
}
=
\frac{ 
\left[
g_z^{\leftarrow}(a^{\theta})
\right]^2
}{4(1-\theta)\ln (\frac{1}{a})}.
$$
By  minimization over $\theta\in (0,1)$, we obtain
$T^*_{c_2}\leq 
{\tilde T}$,
with 
$$
{\tilde T}:=
\inf_{0<\theta<1}\; 
\frac{ 
\left[
g_z^{\leftarrow}(a^{\theta})
\right]^2
}{4(1-\theta)\ln (\frac{1}{a})}.
$$
\vskip0.3cm

It remains to  prove the reverse inequality, namely ${\tilde T}\leq T^*_{c_2}$.  We fix $z\in \R^2$.
We discuss three cases. (The last case may not appear for some functions $g_z$).
\vskip0.3cm

{\em First case}.
Let $\rho\in (0,+\infty)$ such that $a<g_z(\rho) <1$. Then there exists $\theta\in (0,1)$ such that $g_z(\rho)=a^{\theta}$.
Since $g_z^{\leftarrow}(a^{\theta})\leq \rho$   by definition of
$g_z^{\leftarrow}$, we can say that

\begin{equation}\label{tildet}
{\tilde T}
 \leq 
\frac{ 
\left[
g_z^{\leftarrow}(a^{\theta})
\right]^2
}{4(1-\theta)\ln (\frac{1}{a})}
\leq
\frac{ \rho^2}{4\ln (a^{\theta-1})}
=
\frac{ 
\rho^2
}{4 \ln_+ \left(\frac{g_{z}(\rho)}{a}\right)}
=
\frac{ 
\rho^2
}{4 \ln_+ \left(\frac{{\mathcal M}_{z}(\rho)}{L(M)}\right)}
,
\end{equation}
using  the relations 
$a^{\theta-1}=a^{-1} g_z(\rho)=\frac{{\mathcal M}_{z}(\rho)}{L(M)}$.
\vskip0.3cm

{\em Second case}.
Assume that $g_z(\rho) \leq a$, then the inequality \eqref{tildet} holds true 
trivially since the right-hand side  is infinite.
Thus, from the first  and second case, we deduce that
 the inequality \eqref{tildet} holds true  whenever  $g_z(\rho)<1$.
\vskip0.3cm

{\em Third case}.
If it happens that $g_z(\rho)=1$ for some finite $\rho>0$, by the same argument defining $g^{\leftarrow}(m)$ for $m\in (0,1)$, 
we set 
$$
\rho_0=\inf\{\rho^{\prime}>0, g_z(\rho^{\prime})=1\}
={\rm min}\{\rho^{\prime}>0, g_z(\rho^{\prime})=1\}.
$$
So, we have $0< \rho_0\leq \rho$ for $\rho>0$ such that $g_z(\rho)=1$.
For any  sequence $(r_n)$ such that $r_n< \rho_0$ and   $\lim_n r_n=\rho_0$,   then we have  $g_z(r_n)< g_z(\rho_0)=g_z(\rho)=1$.
We can apply the result of the  first  or  second case  just above to  such sequence $(r_n)$, i.e \eqref{tildet}, and  write for all $n$,
$$
{\tilde T}
\leq
\frac{ 
r_n^2
}{4 \ln_+ \left(\frac{{\mathcal M}_{z}(r_n)}{L(M)}\right)}.
$$
Now, taking the limit as $n$ goes to infinity and using  the continuity of the function  $r\mapsto {\mathcal M}_{z}(r)$,
we obtain 
$$
{\tilde T}
\leq
\frac{ 
\rho_0^2
}{4 \ln_+ \left(\frac{{\mathcal M}_{z}(\rho_0)}{L(M)}\right)}
\leq
\frac{ 
\rho^2
}{4 \ln_+ \left(\frac{{\mathcal M}_{z}(\rho)}{L(M)}\right)}.
$$
The last inequality is due to the equality  $g_z(\rho_0)=g_z(\rho)=1$,
i.e. ${\mathcal M}_{z}(\rho)={\mathcal M}_{z}(\rho_0)=M$ and $0< \rho_0\leq \rho$.
\vskip0.3cm

{\em Conclusion}. Gathering the three cases above, we can  assert  that 
 the inequality  \eqref{tildet} holds true for all $\rho>0$. Hence, we obtain 
 ${\tilde T}\leq T_{c_2}^*$  by minimization over $\rho>0$ of  the utmost  right-hand side term of  the inequality \eqref{tildet}.
 Ultimately, we get  ${\tilde T}=T_{c_2}^*$
 since  ${\tilde T}\geq T_{c_2}^*$ has already been proved in the first part of the proof.
 Finally, the last statement   of part  (2)  of Corollary  \ref{corapplica1} holds true by noting that 
$\frac{1}{a}=1+\frac{M-8\pi}{2M}$.
\vskip0.3cm

(3)
 The proof is similar to the proof of the first part  of Corollary   \ref{corapplica1}. 
We assume that $n_0$ has  compact support in $\R^2$, here denoted  by $K$.
Let $T=T^*_{c,z}(n_0)$ as in Theorem \ref{tcalphazero} for fixed $z\in \R^2$.
By  \eqref{heatcondit} of Theorem \ref{tcalphazero}, we  have 
$$
e^{-{i^2_K(z)}/{4T}}.M
=\left(
e^{-\sup_{x\in K}\frac{\vert x -z\vert^2}{4T}}\right).
M
=
\left(
\inf_{x\in K}
e^{-\frac{\vert x -z\vert^2}{4T}}\right).
 \int_{K}  n_0(y)\, dy
 $$
 $$
 =
 \int_{K}   
\inf_{x\in K}
e^{-\frac{\vert x -z\vert^2}{4T}} n_0(y)\, dy
\leq 
  \int_{K} 
  e^{-\frac{\vert y -z\vert^2}{4T}}
  n_0(y)\, dy
$$
$$
=
 \int_{\Ri^2} 
e^{-\frac{\vert y -z\vert^2}{4T}}
n_0(y)\, dy
=
H_{z,n_0}(T)
=
L(M),
$$
with $L(M)=\frac{2M^2}{3M-8\pi}$. This implies that, for all $z\in \R^2$,
$$
T=T^*_{c,z} \leq
\frac{i^2_K(z)}
{4\ln\left(\frac{M}{L(M)}\right)},
$$
due to the fact that  $\ln\left(M/L(M)\right)>0$ since $L(M)< M$ for $M>8\pi$.
By minimization over $z\in \R^2$, the conclusion follows:
$$
T^*\leq T^*_{c} \leq 
T^*_{c_3}:=
\frac{R_0^2}
{4\ln\left(\frac{M}{L(M)}\right)}
=
\frac{R_0^2}{4\ln\left(
1+\frac{M-8\pi}{2M}
\right)}.
$$

An alternative proof  of this inequality \eqref{tctrois} is as follows.
In  the inequality \eqref{tdeuxd},  we can consider  the value  $\theta=0$ (i.e. the limit case  $\theta$ where  goes to $0$). Indeed, 
 $g_z^{\leftarrow}(1)$ is finite for all $z\in \R^2$ (due to the assumption of compact support on $n_0$)
and $g_z^{\leftarrow}(a^{\theta})\leq g_z^{\leftarrow}(1)$
since $t\mapsto g_z^{\leftarrow}(t)$ is non-decreasing and $a^{\theta}< 1$.
So, we can  write
$$  
T^*_{c_2}
\leq
\frac{1}
{
4\ln\left(
1+\frac{M-8\pi}{2M}
\right)
}.
\inf_{z\in \Ri^2,\theta \in (0,1)}
\left(
\frac{ 
\left[
g_z^{\leftarrow}(a^{\theta})
\right]^2
}{(1-\theta)}
\right)
\leq
\frac{1}
{
4\ln\left(
1+\frac{M-8\pi}{2M}
\right)
}.
\inf_{z\in \Ri^2}
\left[
g_z^{\leftarrow}(1)
\right]^2,
$$
Now, it can be easily shown that $g_z^{\leftarrow}(1)=i_K(z)$ for all  $z\in \R^2$. Hence, we get
$$
R_0^2=
\inf_{z\in \Ri^2}
\left[
g_z^{\leftarrow}(1)
\right]^2.
$$
Thus, we   conclude that $T^*_{c_2}\leq T^*_{c_3}$. So, the inequality  \eqref{tctrois}  also follows from  \eqref{tdeuxd}.
This finishes  this second proof.
\vskip0.3cm

The inequality  \eqref{tctroisbis} is a consequence of  Jung's Theorem \cite{wiki}
applied to the two-dimensional case. Indeed, we have 
$R_0\leq D/\sqrt{3}$, where $D$ is the diameter of $K$,
see  Remark \ref{jungth} above.
This concludes the proof of Corollary  \ref{corapplica1}. 
\hfill $\square$
\vskip0.3cm

In the next statement, we  give practical criteria for  obtaining  bounds on $T^*_{c_2}(0)$ of Corollary  \ref{corapplica1}.
Under  natural additional   assumptions on the function $h:=(g_0^{\leftarrow})^2$
appearing in  the inequality   \eqref{tctrois} of Corollary \ref{corapplica1},
we prove the existence of a unique extremum  at $\theta_0\in (0,1)$, or   at $\theta_0=0^+$
of  the infimum over $\theta\in (0,1)$ used in \eqref{tdeuxd}. We also provide an expression of $T^*_{c_2}(0)$ by evaluating
the value of $\theta_0$ by inverting some   function, namely $F$  canonically associated with $g_0^{\leftarrow}$  when $g_0^{\leftarrow}$ is smooth  enough (see below).  
In Section \ref{examples}, we apply   these  estimates  to obtain explicit bounds on the critical time $T^*_{c_2}(0)$ 
for  several  families of examples of  initial data.
 
\begin{cor}\label{corapplica3}
Let $M>8\pi$ and $T^*_{c_2}(0)$ given in \eqref{atheta} of Corollary \ref{corapplica1}, i.e. 
$$ 
T^*_{c_2}(0)
=
\frac{1}
{
4\ln\left(
1/a\right)
}.
\inf_{\theta \in (0,1)}
\frac{ 
h(a^{\theta})
}{(1-\theta)}
,
$$
with $a=\frac{2M}{3M-8\pi}\in (\frac{2}{3}, 1)$ (thus  $\frac{1}{a}=1+\frac{M- 8\pi}{2M}\in (1,\frac{3}{2})$).
Here, the function $h$ is defined by 
$ h(t):=\left[g_0^{\leftarrow}(t)\right]^2$, $t\in (a,1)$. 
Assume that $h$ is continuous and has  derivative $h^{\prime}(t)>0$ for all $t\in (a,1)$.
Let  $F$ be given by
\begin{equation}\label{eqF}
F(X)=X+e^X \left(\frac{h}{h^{\prime}}\right)(e^{-X}),\quad X\in (0, \ln(1/a)).
\end{equation}
We have the next two results.
\begin{enumerate}
\item
If $F$ is non-decreasing and $F(0^+)\geq \ln (1/a)$, i.e.
$ \left(\frac{h}{h^{\prime}}\right)(1^-)\geq \ln(1/a)$, then we have
\begin{equation}\label{tcF1}
T^*_{c_2}(0)=
\frac{h(1^{-})}
{4 \ln(1/a)}.
\end{equation}
\item
If $F$ is continuous and strictly increasing on $ (0, \ln(1/a))$,
and satisfies the next conditions,
\begin{itemize} 
\item[(i)]
$F(0^+) < \ln (1/a) $, i.e. 
$ \left(\frac{h}{h^{\prime}}\right)(1^-) < \ln(1/a)$,
\item[(ii)]
$
F(\left[\ln(1/a)\right]^{-})> \ln(1/a)
$,
i.e. 
$ \left(\frac{h}{h^{\prime}}\right)(a^+)>0$,
\end{itemize} 
then  we have
\begin{equation}\label{tcF2}
T^*_{c_2}(0)=
S\left(\ln \frac{1}{a}\right),
\end{equation}
where 
\begin{equation}\label{funcw}
S(Y):=
\frac{h (e^{-F^{-1}(Y)}
)}
{4 \left[Y-F^{-1}(Y)\right]}
=
\frac{1}{4}
e^{-F^{-1}(Y)}
h^{\prime}(
e^{-F^{-1}(Y)}
)
, \; Y\in {\rm Dom}(F^{-1}),
\end{equation}
and $F^{-1}$ is the inverse function of $F$.  
\end{enumerate}
\end{cor}

We make some comments on these results.

\begin{rem}
The functions $S$ and $h$  are  formally  independent of the mass $M$. These functions 
depends on  the shape of the initial data $n_0$ via the 
repartition function $g_0^{\leftarrow}$. 
Both  expressions   \eqref{tcF1} and \eqref{tcF2} are expressed respectively in terms of the  functions $h$ and $S$  of the variable $Y$
evaluated at $Y_0=Y_0(M):=\ln(\frac{1}{a})=\ln\left(1+\frac{M-8\pi}{2M}\right)$ where $M$  is  the mass of the initial data $n_0$ with $M>8\pi$. 
Note that we  have $\ln(\frac{1}{a})\in (0, \ln(3/2))$ where  $\ln(3/2)\sim 0.40$, thus $\ln(\frac{1}{a})\in (0,0.41)$ for all $M>8\pi$. 
These comments are provided so that the reader can get an idea of the quantities involved in the computations.
\end{rem}

\begin{rem}
In Section \ref{examples}, we evaluate  the function $F$ for some families  of  examples of initial data.  
But unfortunately, the inverse function  $F^{-1}$ can not  always  be  given  explicitly  as  simple functions. 
Nevertheless for some examples, the  function  $F^{-1}$  can be bounded above and below by explicit functions useful for estimating $T^*_c(n_0)$.
Hence, it is also useful  to provide  an upper bound on $T^*$.
Since only the evaluation of $F^{-1}$ at   $Y_0=\ln( \frac{1}{a})$ is involved in \eqref{tcF2},
 the quantity $F^{-1}(\ln  \frac{1}{a})$ may certainly  be estimated by numerical methods 
 as  the unique zero of the function  $Q:=F+\ln a$. 
 But  we shall not continue in that direction  in this paper.
\end{rem}

\begin{rem}
 The map $M\in (8\pi ,+\infty)\rightarrow (\frac{2}{3}, 1)$  defined by 
$a=a(M)=\frac{2M}{3M-8\pi}$ is a decreasing  bijection. Hence $h$ is defined at least on $(a,1)$ 
for all $a\in (\frac{2}{3},1)$,  i.e. on $ (\frac{2}{3},1)$ itself.
This  implies that $F$ defined in Corollary \ref{corapplica3} must have a domain containing $(0,\ln (\frac{3}{2}))$.
So,  the  function $F$ does not depend on $a$. In fact, the function $F$ is often defined on the whole interval  $(0,+\infty)$.
\end{rem}

\begin{rem}
The assumption $ \left(\frac{h}{h^{\prime}}\right)(a^+)>0$ in the second part of the corollary can be replaced by 
the  simple assumption $h^{\prime}(a^+)>0$. Indeed, it is enough to show that we always have  
$h(a^+)>0$.  This is deduced as follows.
Since we have $h=(g_0^{\leftarrow})^2\geq 0$,   we have $h(a^+)\geq 0$.
 It remains to show that $h(a^+)\neq 0$.  To prove this, we suppose the contrary, i.e. $h(a^+)= 0$. 
 Thus, for any sequence $(a_n)$ such that $a\leq a_n<1$ and 
 $\lim_n a_n=a$, then we have  $h(a^+)= \lim_n h( a_n)=0$. This implies that $ \lim_ng_0^{\leftarrow}( a_n)=0$.
Thus, on one hand  we get 
$$
\lim_n g_0 ( g_0^{\leftarrow}(a_n))= g_0(0)=0,
$$
by  continuity of $g_0$.
On the other hand, we have
$$
\lim_n \,(g_0 \circ g_0^{\leftarrow})(a_n)=\lim_n  a_n=a,
$$
due to  the fact that $g_0^{\leftarrow}$ is a right-inverse of $g_0$.
So, we deduce that  $a=\frac{2M}{3M-8\pi}=0$, which contradicts the assumption $M>8\pi$. 
\end{rem}

\begin{rem}
For  any $z\in \R^2$, and under similar assumptions on $h_z:=(g_z^{\leftarrow})^2$  as those  imposed on $h=h_0$,  
we can prove    similar  estimates  as \eqref{tcF1} and  \eqref{tcF2}  for $T^{*}_{c_2}(z)$ 
with the corresponding $F_z$.
\end{rem}

{\bf Proof of Corollary \ref{corapplica3} }.
We first  need some preparation.
Let $a$ and $h:=h_0$  defined as  in Corollary \ref{corapplica3} ($z=0\in \R^2$). 
We assume that $h$ is continuous and  it has a positive  derivative.
We set $V(\theta):=\frac{ h(a^{\theta})}{(1-\theta)}$ with $\theta\in (0,1)$.
The derivative of the function $V$ satisfies the  following  equation 
$$
(1-\theta)^2V^{\prime}(\theta)=
(1-\theta)(\ln a) h^{\prime}(a^{\theta})a^{\theta}+h(a^{\theta}).
$$
We make the following change of variables $X=X(\theta)=-\theta \ln a=\theta \ln( 1/ a)$. 
The map $X$ is an increasing bijection from $(0,1)$ onto
$(0,\ln (1/a))$. We deduce that $a^{\theta}=e^{-X}$, and
$$
(1-\theta)^2V^{\prime}(\theta)=
(1-\theta)(\ln a) h^{\prime}(e^{-X})e^{-X}+h(e^{-X})
$$
$$
=
(\ln a) h^{\prime}(e^{-X})e^{-X}+ Xh^{\prime}(e^{-X})e^{-X}
+ h(e^{-X})
=
h^{\prime}(e^{-X})e^{-X}
Q(X) 
,
$$
with $Q(X):=F(X)+\ln a$, and $F$ given by \eqref{eqF}.
From the assumption $h^{\prime}>0$, we get
$Q(X)>0 \iff V^{\prime}(\theta)>0$  and $Q(X)=0 \iff V^{\prime}(\theta)=0$,
where the relationship  between $X$ and $\theta$ is given by $X=-\theta\ln a$ . 
Now we are in position to prove our statements.
\vskip0.3cm

(1) 
Under the assumptions of (1) of  Corollary \ref{corapplica3}, we have 
$Q(0^+)=F(0^+)+\ln a\geq 0$, and  $Q$ (also  $F$) is non-decreasing.
This implies that 
$Q(X)\geq Q(0^+)\geq 0$, for all $X \in (0,\ln (1/a))$.
Thus,  we have $V^{\prime}(\theta)\geq 0$ for all $\theta\in (0,1)$, i.e.
$V$ is  non-decreasing. Hence, we deduce  that
$$
\inf_{0<\theta<1} V(\theta)=V(0^+)=h(1^{-}).
$$
Then the formula  \eqref{tcF1} follows from  \eqref{tdeuxd} of Corollary \ref{corapplica1} with $z=0$.
This proves the first statement.
\vskip0.3cm

(2)\, Under the assumptions of (2) of  Corollary \ref{corapplica3},   the function $Q(X)=F(X)+ \ln a $  is a   strictly increasing continuous function,
and by (i) and (ii), we have successively  $Q(0^+)<0$ and $Q(\left[\ln(1/a)\right]^{-}) >0$.
By  the intermediate value  theorem,
there exists a unique  zero 
$X_0\in (0, \ln (1/a))$ of $Q(X)=F(X)+\ln a$, i.e. $F(X_0)=\ln  (1/a)$, and finally $X_0= F^{-1}(\ln \frac{1}{a})$.
Indeed,  $F$ is  bijective from $ (0, \ln (1/a))$ onto its range
which contains  $\ln (1/a)$ since $F(0^+) < \ln \frac{1}{a} < F(\left[\ln(1/a)\right]^{-})$ by assumptions (i) and (ii).
We also  obtain  
$$
 Q(X)<  0
\iff
0<X < X_0,
\quad
{\rm and}\quad Q(X)= 0
\iff
X =X_0,
$$
which is equivalent to 
 $$
 V^{\prime}(\theta)<0\iff 0<\theta<\theta_0,
 \quad
 {\rm and}\quad  V^{\prime}(\theta)= 0
\iff
\theta =\theta_0,
 $$
with the relations  $X=\theta \ln(1/a)$ and $X_0=\theta_0\ln(1/a)$.
 This implies that  the infimum of $V(\theta)$ over $\theta\in (0, 1)$ is attained at this point $\theta_0\in (0, 1)$, i.e.
 $$
 \inf_{0<\theta<1} V(\theta)= V(\theta_0),
  $$
 with 
 $$
 \theta_0=
 \frac{X_0}{\ln (1/a)}
 =
  \frac{F^{-1}(\ln  \frac{1}{a})}
{\ln  (1/a)}.
   $$
Thus, we derive  the  next formula 
$$ 
T^*_{c_2}(0)
=
\frac{1}
{
4\ln\left(
1/a\right)
}.
V(\theta_0)
=
\frac{1}
{
4\ln\left(
1/a\right)
}.
\frac{ 
h(a^{\theta_0})
}{(1-\theta_0)}
=
\frac{h\left(
e^{-F^{-1}(\ln \frac{1}{a})}
\right)}
{4 \left[\ln(\frac{1}{a})-F^{-1}(\ln \frac{1}{a})\right]}.
$$
Now let  $X=F^{-1}(Y)$ with $Y\in {\rm Dom}(F^{-1}):={\rm Im} \, F$ in the definition \eqref{eqF} of $F(X)$.
We get the following  expression
$$
Y-F^{-1}(Y)= e^{F^{-1}(Y)}\left( \frac{h}{h^{\prime}}\right)  (e^{-F^{-1}(Y)}).
$$
From which  we deduce both formulas for the function $S$, namely
$$
S(Y):=\frac{h(e^{-F^{-1}(Y)})}
{4 \left[Y-F^{-1}(Y)\right]}
=
\frac{1}{4} e^{-F^{-1}(Y)}
h^{\prime}(e^{-F^{-1}(Y)}).
$$ 
Finally, this yields the expected result
$$
T^*_{c_2}(0)=
S\left(\ln \frac{1}{a}\right)
$$
where  $S(Y)$ is described just above, and $Y=\ln(\frac{1}{a})$.
Statement (2) is proved and the proof of Corollary \ref{corapplica3} is now complete.
\hfill$\square$
\vskip0.3cm

Now, we present  another consequence of Theorem  \ref{tcalphazero} under an additional assumption of finite centralized and normalized 
$\beta$-variance of the initial data $n_0$ defined as follows. We assume that $n_0$ has a moment of order $\beta\geq 1$, i.e.
$$
\int_{\Ri^2} \vert x\vert^{\beta} n_0(x)\, dx <+\infty.
$$
\noindent
Then, we denote the (normalized) barycenter of $n_0$ by $B_0=\frac{1}{M} \int_{\Ri^2}  x. n_0(x)\,dx$. This is a well defined  vector in $\R^2$, 
since by H\"older's inequality, we obtain
$$
\frac{1}{M} \int_{\Ri^2} \vert  x \vert  n_0(x)\,dx
\leq
\left(\frac{1}{M} \int_{\Ri^2} \vert  x \vert^{\beta}  n_0(x)\,dx\right)^{1/\beta} <+\infty,
$$
for any $1\leq \beta<+\infty$.
We also denote by  $V_{\beta}(n_0)$  the centralized and normalized $\beta$-variance of $n_0$ defined by 
$$
V_{\beta}(n_0):=
\left[
\frac{1}{M}
\int_{\Ri^2} \vert x-B_0\vert^{\beta} n_0(x) \,dx\right]^{2/\beta}.
$$
Let us mention some basic properties  of the  $\beta$-variance of $n_0$   useful for this paper.
\vskip0.3cm

(a)
The $\beta$-variance, seen as a function of $\beta$, i.e.  $\beta\in [2,+\infty)\mapsto V_{\beta}(n_0)$, is non-decreasing. 
Indeed, by  H\"older's  inequality,  we have  for  all $2\leq \beta\leq \gamma$,
 $$
 V_{\beta}(n_0)
 =\left[
\frac{1}{M}
\int_{\Ri^2} \vert x-B_0\vert^{\beta} n_0(x) \,dx\right]^{2/\beta}
\leq
\left[
\frac{1}{M}
\int_{\Ri^2} \vert x-B_0\vert^{p\beta} n_0(x) \,dx\right]^{2/p \beta}
=  V_{\gamma}(n_0),
$$
where    $p=\gamma/\beta\geq 1$. 
As a  particular case, we get
$$
V_{2}(n_0)=
\frac{1}{M}
\int_{\Ri^2} \vert x-B_0\vert^2 n_0(x) \,dx
\leq
V_{\beta}(n_0)=
\left[
\frac{1}{M}
\int_{\Ri^2} \vert x-B_0\vert^{\beta} n_0(x) \,dx\right]^{2/\beta},
$$
for  all  $2\leq \beta<+\infty$.
As a consequence, if $V_{\beta}(n_0)$  is finite for some $\beta\geq 2$ then $V_2(n_0)$  is also finite.
\vskip0.3cm

(b) 
We have the following  infimum estimate
\begin{equation}\label{lowbeta}
\frac{1}{4}
V_{\beta}(n_0)
\leq 
\inf_{z\in \Ri^2}
\left[\frac{1}{M}
\int_{\Ri^2} \vert x-z\vert^{\beta} n_0(x) \,dx
\right]^{2/\beta}
\leq
V_{\beta}(n_0).
\end{equation}
\vskip0.3cm

(c)
In  case $\beta=2$,
 we have an equality of the infimum with the upper bound in \eqref{lowbeta},
\begin{equation}\label{lowdeux}
\inf_{z\in \Ri^2}
\left[\int_{\Ri^2} \vert x-z\vert^{2}\,\frac{n_0(x)}{M}\,dx
\right]
=
V_{2}(n_0)
=
\int_{\Ri^2} \vert x- B_0\vert^{2}\,\frac{n_0(x)}{M}\,dx
.
\end{equation}
Note that the quantity $V_{2}(n_0)$ is simply  denoted by $V(0)$ in  the first section of this paper.
\vskip0.3cm

The proof of the lower bound of \eqref{lowbeta}
for $\beta\geq 2$ is as follows. (The upper bound is obvious taking $z=B_0$).
From (discrete)  H\"older's inequality, we get
$$
(a+b)^{\beta}
\leq 
2^{\beta-1}
\left(
a^{\beta}
+
b^{\beta}\right),\quad a,b\geq 0.
$$
Then, from this inequality and triangle inequality,  
we deduce that
$$
\vert x-B_0\vert^{\beta}\leq
2^{\beta-1}
\left(
\vert x-z\vert^{\beta}
+
\vert z-B_0\vert^{\beta}\right),
$$
for all  $x,z\in \R^2$. We multiply this expression by the non-negative function $n_0/{M}$  and  integrate it  with respect to $x$ over $\R^2$,  
we then obtain
$$
\int_{\Ri^2} \vert x-B_0\vert^{\beta}\,\frac{n_0(x)}{M}\,dx
\leq
2^{\beta-1}
\left(
\int_{\Ri^2} \vert x-z\vert^{\beta}\,\frac{n_0(x)}{M}\,dx
+
\vert z-B_0\vert^{\beta}\right).
$$
On the other hand, we have 
$$
\vert z-B_0\vert^{\beta}
=\vert \int_{\Ri^2} (z-x)
\,\frac{n_0(x)}{M}\,dx\vert^{\beta}
\leq
\left( \int_{\Ri^2} \vert z-x\vert
\,\frac{n_0(x)}{M}\,dx\right)^{\beta}
\leq
 \int_{\Ri^2} \vert z-x\vert^{\beta}
\,\frac{n_0(x)}{M}\,dx.
$$
The last inequality is obtained by applying Jensen's inequality to the probability measure $\frac{n_0(x)}{M}\,dx$ (or  simply by applying  H\"older's inequality).
It follows that
$$
\int_{\Ri^2} \vert x-B_0\vert^{\beta}\,\frac{n_0(x)}{M}\,dx
\leq
2^{\beta}
\int_{\Ri^2} \vert x-z\vert^{\beta}\,\frac{n_0(x)}{M}\,dx,
$$ for all $z\in \R^2$. 
We conclude the proof of the lower bound of \eqref{lowbeta}
 by raising the expression just   above  to the power  of 
$2/\beta$ and by taking the infimum over $z\in \R^2$.
\vskip0.3cm

The next result   emphasizes  the importance of the role of the  normalized barycenter $B_0$  of  $n_0$
in the issue of estimating $T^*_c(n_0)$ (hence $T^*$ also)  under a $\beta$-moment  assumption on the initial data $n_0$.

\begin{cor}\label{corapplica2} (Finite Variance)
\begin{enumerate}
\item
Assume that  the initial data  $n_0$ has a 2-moment. Then,  we have the following estimate
\begin{equation}\label{tstarparti2}
T^*\leq T^*_c(n_0)\leq
T_{c_4}^*:=
\frac{
V_{2}(n_0) 
}
{4\ln\left(1+\frac{M-8\pi}{2M}\right)},
\end{equation}
where $V_{2}(n_0)$ is the normalized variance of $n_0$.
 ($V_{2}(n_0)$ is also denoted by $V(0)$  in \eqref{bupbdp} above).
\vskip0.3cm

\item
More generally, assume that $n_0$ has a $\beta$-moment with $\beta \geq 2$. Then we have
\begin{equation}\label{tstarparti3}
T^*\leq T^*_c(n_0)\leq
T_{c_5}^*:=
\frac{
\inf_{z\in \Ri^2}
\left[\int_{\Ri^2} \vert x-z\vert^{\beta}\,n_0(x)\,dx
\right]^{2/\beta}
}
{4M^{2/\beta}\ln\left(1+\frac{M-8\pi}{2M}\right)}.
\end{equation}
In particular,  we obtain
\begin{equation}\label{tstarparti}
T^*\leq  T^*_c(n_0)\leq
\frac{
V_{\beta}(n_0)
}
{4\ln\left(1+\frac{M-8\pi}{2M}\right)}, \quad \beta\geq 2.
\end{equation}
\end{enumerate}
\end{cor}

We make some comments on these results.
\begin{rem}
The estimate \eqref{tstarparti3}  is apparently better than the estimate \eqref{tstarparti} (except  the case $\beta=2$ for which we have equality).
See  Remarks \eqref{lowbeta} and  \eqref{lowdeux}.
\end{rem}

\begin{rem}
Here again,  the  term
$\left[\ln(1/a)\right]^{-1}
=
\left[
\ln\left(1+\frac{M-8\pi}{2M}\right)
\right]^{-1}
$ 
appears in our estimates of $T^*_c(n_0)$ as in Corollary  \ref{corapplica1} and  Corollary  \ref{corapplica3}. 
Note that this term tends to $+\infty$ when $M\rightarrow 8\pi^+$, see comments in first section .
\end{rem}

\begin{rem}
In  the case of  finite second-moment  for the initial data $n_0$,
we can  compare two asymptotic results.
From  \eqref{bupbdp2}, we first recall that we can write  the next upper bound 

\begin{equation}\label{virialcomp}
T^*\leq \frac{2\pi V_2(n_0)}{M-8\pi},
\end{equation}
 for all $M>8\pi$.
On the other hand, from \eqref{tstarparti2} we have 
\begin{equation}\label{compar1}
T^*\leq T_c^*(n_0)\leq 
\frac{
V_{2}(n_0)
}
{4\ln\left(1+\frac{M-8\pi}{2M}\right)}, \quad M>8\pi.
\end{equation}

\begin{enumerate}
\item[(i)]
{\em When $M$ is closed to $8\pi$}.
The right-hand side  of \eqref{compar1}  is estimated  formally as follows
$$
\frac{
V_{2}(n_0)
}
{4\ln\left(1+\frac{M-8\pi}{2M}\right)}
\sim
\frac{4\pi V_{2}(n_0)}
{M-8\pi}, \quad  M\rightarrow 8\pi^+.
$$
(i.e. when considering $V_{2}(n_0)$ is fixed).
So, the estimate of  $T^*$ obtained from   \eqref{compar1} is twice larger than  the  estimate 
\eqref{virialcomp} as $M$ tends to $8\pi^+$.
\vskip0.3cm
\item[(ii)]
{\em The asymptotic case $M\rightarrow +\infty$}.
From \eqref{compar1}, we have formally
$$
\frac{
V_{2}(n_0)
}
{4\ln\left(1+\frac{M-8\pi}{2M}\right)}
\sim
\frac{V_{2}(n_0)}
{
4\ln(3/2)},
$$
as $M$ tends to infinity, which is clearly  not as good  as  the estimate  \eqref{virialcomp} of $T^*$. 
Indeed,
 the  right-hand side term of   \eqref{virialcomp} tends to zero as $M$ tends to infinity.
\vskip0.3cm
The fact that the estimate \eqref{virialcomp} of $T^*$ is  better than \eqref{compar1}  is due to   a direct  study 
of the evolution  in time of the second moment  of the solution $(n_t)$. 
For  this  particular case, the general  approach used in this paper with the heat kernel is not  specific 
enough for obtaining  an accurate result.
\end{enumerate}
\end{rem}

\begin{rem}
Assume that $n_0$ has compact support. 
Let $K$ denote  the support of $n_0$. Then for any $R_0>0$, and any  $z_0\in \R^2$
such that $K \subset B^{\prime}(z_0, R_0)$ (closed ball).
Then,  it is easy to show  that
$$\inf_{z\in \Ri^2}
\left[\frac{1}{M}
\int_{\Ri^2} \vert x-z\vert^{\beta} n_0(x) \,dx
\right]^{2/\beta}
\leq
\left[\frac{1}{M}
\int_{B^{\prime}(z_0, R_0)} \vert x-z_0\vert^{\beta} n_0(x) \,dx
\right]^{2/\beta}
\leq R_0^2,
$$
for all $\beta \geq 2$.
Thus,  the inequality \eqref{tstarparti3}   immediately  implies  the inequality \eqref{tctrois} of Corollary \ref{corapplica1}.  
Note also  that  the $\beta$-variance 
$V_{\beta}(n_0)$ is finite since we have
$\frac{1}{4}V_{\beta}(n_0)\leq R_0^2$  obtained from  the inequality \eqref{lowbeta}. 
More precisely for the case $\beta=2$, we  obtain the following family of inequalities,
$$
V_{2}(n_0)
=\frac{1}{M}
\int_{\Ri^2} \vert x\vert^2 n_0(x) \,dx
-\vert B_0\vert^2
\leq R_0^2-\vert B_0\vert^2
\leq  R_0^2,
$$
where $B_0=\frac{1}{M} \int_{\Ri^2} x. n_0(x) \,dx$ ($\in \R^2$)  is the (normalized) barycenter of $n_0$. 
From these last remarks, we can see  that the bounds on $T^*_c(n_0)$  obtained with the variance are stronger in general 
than the one obtained with the minimal radius $R_0$ for the  case where $n_0$ has compact support. 
Note also that the scope of applications  with the variance is wider. Indeed, the variance of $n_0$ 
can be finite without  compact support  condition for $n_0$.
See  the examples of   Section \ref{examples}.
\end{rem}

\noindent{\bf Proof of Corollary \ref{corapplica2}} 

Assume that  $n_0\in L^1$ with $n_0\neq 0$ and it has a $\beta$-moment with $\beta\geq 2$.
 So, the $\beta$-variance $V_{\beta}(n_0)$ is finite.
To prove the corollary, it is enough to prove  the inequality  \eqref{tstarparti3}.  Indeed, this  implies  \eqref{tstarparti2}  with  $\beta=2$ by using
the relation \eqref{lowdeux}. This also gives \eqref{tstarparti} by  using the upper bound of 
\eqref{lowbeta} proved independently and above the statement of Corollary \ref{corapplica2}.
 The lower bound in \eqref{lowbeta} shows the sharpness of the result 
 \eqref{tstarparti} with respect to   \eqref{tstarparti3} (up to the multiplicative constant $1/4$).
\vskip0.3cm

Let $M>8\pi$, $z\in \R^2$, and   $\beta \geq 2$ be fixed. 
By definition of $T:=T^*_c(n_0)$ as defined   in \eqref{tcnormal}  of Theorem \ref{tcalphazero},
 we have  $T\leq T^*_{c,z}$ for all $z\in \R^2$.
Since  $s\mapsto H_{z,n_0}(s)$ is non-decreasing, then we deduce that
$$
\int_{\Ri^2} \exp\left(-\frac{\vert x -z\vert^2}{4T}\right)n_0(x) \, dx
=
H_{z,n_0}(T)\leq 
H_{z,n_0}(T^*_{c,z})=\frac{2M^2}{3M-8\pi}.
$$
We can rewrite this inequality   in the following form,
$$
\int_{\Ri^2}\Psi\left(\left[\frac{\vert x -z\vert}{2\sqrt{T}}\right]^{\beta}\right)  d\mu(x)
=
\int_{\Ri^2} \exp\left(-\frac{\vert x -z\vert^2}{4T}\right)\frac{n_0(x)}{M} \, dx\leq
\frac{2M}{3M-8\pi},
$$
with  $\Psi(r)=\exp(-r^{\gamma})$, $\gamma=2/\beta \leq 1$ and $d\mu(x)=\frac{n_0(x)}{M}\, dx$. 
Now, because  the function $\Psi$ is convex and $\mu$ is a probability measure,  
we can  apply  Jensen's inequality   and  deduce that
$$
\Psi\left( \int_{\Ri^2}\left[\frac{\vert x -z\vert}{2\sqrt{T}}\right]^{\beta} d\mu(x)\right) 
\leq
\int_{\Ri^2}\Psi\left(\left[\frac{\vert x -z\vert}{2\sqrt{T}}\right]^{\beta}\right)  d\mu(x).
$$
Thus, we obtain
\begin{equation}\label{bothsid}
\Psi\left( \int_{\Ri^2}\left[\frac{\vert x -z\vert}{2\sqrt{T}}\right]^{\beta} d\mu(x)\right) 
\leq
a,
\end{equation}
where   $a=\frac{2M}{3M-8\pi}< 1$ (since  $M>8\pi$).
The function $\Psi$ is strictly decreasing from $(0,+\infty)$ onto $(0,1)$, 
then it is invertible with $\Psi^{-1}$  its decreasing inverse.
By applying $\Psi^{-1}$ on both sides of the inequality \eqref{bothsid} (and reversing the inequality), 
this  leads to
$$
\frac{1}{2^\beta T^{\beta/2} }
\int_{\Ri^2} \vert x -z\vert^{\beta}\frac{n_0(x)}{M}\, dx
=
\int_{\Ri^2}\left[\frac{\vert x -z\vert}{2\sqrt{T}}\right]^{\beta} d\mu(x)
 \geq 
 \Psi^{-1}(a)>0.
 $$
 This yields,
$$
T^*_c(n_0)=:T
\leq 
\frac{1}{ 4 [\Psi^{-1}(a)]^{2/\beta}}
\left[
\int_{\Ri^2} \vert x -z\vert^{\beta}\frac{n_0(x)}{M}\, dx
\right]^{2/\beta}.
 $$
It is easy to verify  that $\left[\Psi^{-1}(a)\right]^{2/\beta}=\ln(1/a)$.
For all   $z\in \R^2$, we finally get 
$$
T^*_c
\leq 
\frac{1}{ 4 \ln (1/a)}
\left[
\int_{\Ri^2} \vert x -z\vert^{\beta}\frac{n_0(x)}{M}\, dx
\right]^{2/\beta}.
$$
This implies the inequality \eqref{tstarparti3} after minimizing over  $z\in \R^2$.
This concludes the proof of Corollary \ref{corapplica2}.
\hfill $\square$
\vskip0.3cm
 
\begin{rem}
With the same  argument of proof as in  the proof of  Corollary \ref{corapplica2},  we can also show that
$$
T^*_{c,z}
\leq 
\frac{1}{ 4 \ln (1/a)}
\left[
\int_{\Ri^2} \vert x -z\vert^{\beta}\frac{n_0(x)}{M}\, dx
\right]^{2/\beta},
$$
for all $z\in \R^2$.  Obviously, this also implies the  inequality \eqref{tstarparti3}.
\end{rem}
\vskip0.3cm

\begin{rem}
A variant of the  proof  of \eqref{tstarparti}, for any $\beta\geq 2$,  can be  given  
using the inequality  \eqref{tstarparti2}  proved as   the inequality  \eqref{tstarparti}  for the particular case  $\beta=2$, 
 and  using the fact that the map
$$
\beta\in (0,+\infty) \rightarrow {\tilde V}_{z}(\beta):=
\left( \int_{\Ri^2} \vert x -z\vert^{\beta}\frac{n_0(x)}{M}\, dx\right)^{2/\beta}
$$
 is non-decreasing for all  $z\in \R^2$.
Indeed, this last fact can be shown as follows.
For any $0<\gamma \leq \beta=p\gamma$, we have by Jensen's inequality  (or H\"older's inequality) the next inequality 
$$
\left(\int_{\Ri^2} \vert x -z\vert^{\gamma}\frac{n_0(x)}{M}\, dx\right)^p
\leq 
\int_{\Ri^2} \vert x -z\vert^{\gamma p}\frac{n_0(x)}{M}\, dx= \int_{\Ri^2} \vert x -z\vert^{\beta}\frac{n_0(x)}{M}\, dx,
$$
since $p=\frac{\beta}{\gamma}\geq 1$. So, 
we immediately deduce the inequality 
 ${\tilde V}_{z}(\gamma)\leq {\tilde V}_{z}(\beta)$ 
 for all  $z\in \R^2$.
Thus, the inequality  \eqref{tstarparti}   follows from  \eqref{tstarparti2} for all  $\beta\geq \gamma=2$ by taking $z=B_0$,
since $V_{2}(n_0)={\tilde V}_{B_0}(2)\leq {\tilde V}_{B_0}(\beta)=V_{\beta}(n_0)$. 
\end{rem}

\section{Lower bounds on $T^*_c(n_0)$}\label{lowerboundsect}
The aim of this section is to provide a lower bound  on the critical time 
$T^*_c(n_0)$ (but  unfortunately  not on $T^*)$ under some additional $L^p$ properties of the initial data $n_0$.
Note that during  the proof, we have the opportunity to use the best bounds of the  $L^p-L^q$  embedding theorem of the heat semigroup.
For the definition of  $T^*_c(n_0)$, see  \eqref{tcnormal} above.

\begin{theo}\label{lowe1}
Assume that $n_0\in L^1\cap L^p$ for some $1<p\leq +\infty$
and $M=\vert \vert n_0\vert\vert_1>8\pi$.
The following lower bounds on  $T^*_c(n_0)$ hold true.
 \begin{enumerate}
\item
We have
\begin{equation}\label{lowerbd}
T^*_{c_5}
:=
\frac{1}{4\pi}
\sup_{1<q\leq p} q'
\left[ \frac{L(M)}{  \vert \vert n_0\vert\vert_q}
\right]^{q'}
\leq 
T^*_c(n_0),
\end{equation}
where $L(M)=\frac{2M^2}{3M-8\pi}$ 
and 
$\frac{1}{q}+\frac{1}{q^{\,\prime}}=1$.
\item
Let $p_0=\left[ \ln (\frac{2e}{3})\right]^{-1}$, ($p_0 \sim 1,682$).
 \begin{enumerate}
\item
If $p_0\leq p\leq +\infty$ or  if $1<p<p_0$ and $M\leq \frac{8\pi}{3 -2e^{(1-1/p)}}$, then we have
\begin{equation}\label{lowerbd2}
\frac{(\pi e)^{-1}}
{4 \ln\left(1+\frac{M-8\pi}{2M}\right)}
\left[ \frac{M}{  \vert \vert n_0\vert\vert_p}
\right]^{p'}
\leq 
T^*_c(n_0).
\end{equation}
In particular for  $p=+\infty$, we have 
\begin{equation}\label{lowerbd2bis}
\frac{(\pi e)^{-1}}
{4 \ln\left(1+\frac{M-8\pi}{2M}\right)}
\left[ \frac{M}{  \vert \vert n_0\vert\vert_{\infty}}
\right]
\leq 
T^*_c(n_0).
\end{equation}

\item
If  $1<p<p_0$ and $M> \frac{8\pi}{3 -2e^{(1-1/p)}}$, 
then we have
\begin{equation}\label{lowerbd3}
\frac{p'}
{4 \pi}
\left[ \frac{L(M)}{  \vert \vert n_0\vert\vert_p}
\right]^{p'}
\leq 
T^*_c(n_0),
\end{equation}
where $\frac{1}{p}+\frac{1}{p^{\,\prime}}=1$.
\end{enumerate}
\end{enumerate}
\end{theo}

For  applications of Theorem \ref{lowe1}, see Section \ref{examples}.
\vskip0.3cm

We make some comments on these results.

\begin{rem}
If $n_0\in L^1\cap L^p$ then $n_0\in L^q$ for all $q\in [1,p]$   by interpolation (generalized H\"older's inequality).
Thus, the expression $T^*_{c_5}$ of the left-hand side of  \eqref{lowerbd}
is well-defined.
\end{rem}

\begin{rem}
 By $L^{\infty}-L^{\infty}$ contraction of the heat semigroup $(e^{t\Delta})_{t\geq 0}$,   
 the following sharp inequality in $L^{\infty}$ is well-known, 
$$
 \vert\vert p_s\star f \vert\vert_{\infty}\leq \vert\vert f\vert\vert_{\infty},
 $$
 for all $f\in L^{\infty}$ and  all  $s>0$.
 Let $T=T^*_{c,z}(n_0)$ with $z\in \R^2$ be fixed as in Theorem \ref{tcalphazero} (2),
  and $n_0\in L^1\cap L^{\infty}$.
As a consequence of the above inequality, we can deduce that  
$$
L(M)=H_{z,n_0}(T)=(4\pi T)\,  p_T\star n_0 (z)
\leq  4\pi T \vert\vert n_0\vert\vert_{\infty},
$$
where $L(M)=\frac{2M^2}{3M-8\pi}$  and $M=\vert\vert  n_0\vert\vert_1$.
Thus, we obtain 
$$
\frac{L(M)}{4\pi  \vert\vert n_0 \vert\vert_{\infty}}
\leq T^*_{c,z}(n_0),
$$
 for all $z\in \R^2$.
Hence, by minimization over $z\in \R^2$ and replacing  $L(M)$ by its explicit  expression,  we can write
\begin{equation}\label{lowerbd4}
\frac{\pi^{-1}}
{4 \left(1+\frac{M-8\pi}{2M}\right)}.
\left[\frac{M}{ \vert\vert n_0\vert\vert_{\infty}}\right]
\leq 
\inf_{z\in \Ri^2} T^*_{c,z}(n_0)=T^*_c(n_0),
\end{equation}
for all $n_0\in L^1\cap L^{\infty}$. 
\vskip0.3cm

We can prove that the inequality \eqref{lowerbd2bis}
 is  stronger than \eqref{lowerbd4}.
Indeed,  we  can easily prove  the  next sharp pointwise inequality,
$$
\frac{1} {1+u} \leq
{\kappa} \frac{e^{-1}}{\ln (1+u)}
\leq
 \frac{e^{-1}}{\ln (1+u)}
, \quad 0\leq u\leq \frac{1}{2},
$$
where $\kappa=(2e/3)\ln(3/2)\sim 0,735<1$.  But note that a reverse inequality doesn't hold.
 Indeed, there exists no $\delta>0$ such that
$$
 \frac{{\delta} }{\ln (1+u)}
\leq
\frac{1} {1+u},
 \quad 0< u\leq \frac{1}{2}.
$$
(Just by taking $u\rightarrow 0^+$).
Then, we can apply the  inequalities just above with $u=\frac{M-8\pi}{2M}$ with $M>8\pi$, since $0<u<1/2$,
 which shows that the inequality  \eqref{lowerbd2bis}  is  stronger  than the inequality  \eqref{lowerbd4}.
The reason  of  this improvement   \eqref{lowerbd2bis} with  respect to \eqref{lowerbd4} is that 
$n_0\in L^1\cap L^{\infty}$ implies  $n_0\in L^p$ for any $1\leq p\leq {\infty}$ 
by interpolation, not only $n_0\in L^{\infty}$. This allows us  to use all the scale of $L^p$-spaces, and not only $L^1$ and $L^{\infty}$.
For details, see the proof of Theorem \ref{lowe1} below.
\end{rem}

\begin{rem}
 As a consequence of  \eqref{lowerbd2bis}, the singularity of the ratio $\frac{1}
{4 \ln\left(1+\frac{M-8\pi}{2M}\right)}$ as $M$ tends to $ 8\pi^+$ leads to the following  expected consequence,
$$
\lim_{M\rightarrow 8\pi^+} T^*_c(n_0)=+\infty,
$$ for all $n_0\in L^1\cap L^{p}$,
$1<p\leq +\infty$.
Note that this result can not be obtained directly  from \eqref{lowerbd4}.
\end{rem}

\begin{rem}
A lower bound on $T^*$ of  similar form as in \eqref{lowerbd} for   $T^*_c$ has been obtained in \cite[Th. 1 (i)]{koz-sug}
when $n_0\in L^p$ with  the restriction $1<p<2$. See also  \cite[Prop.1.1]{koz-sug} for $n_0\in L^{\infty}$. 
\end{rem}

{\bf Proof of Theorem \ref{lowe1}}. (1)
We first recall a general result about $L^p-L^q$ contractions  of the heat semigroup $(e^{t\Delta})_{t\geq 0}$ 
defined by convolution with the heat kernel $p_t$. On $\R^n$, we have  the next inequality,
$$
\vert\vert p_t\star  f
\vert\vert_q=\vert\vert e^{t\Delta} f
\vert\vert_q\leq
C(n,p,q) t^{-\frac{n}{2}(\frac{1}{p}-\frac{1}{q})} \vert\vert f\vert\vert_p,  \quad t>0, \; f\in L^p(\R^n),
$$
 with  sharp constants 
 $$
 C(n,p,q)
 =\left( \frac{C_p}{C_q}\right) ^n
 \left[ \frac{4\pi }{(\frac{1}{p}-\frac{1}{q})}\right]^{-\frac{n}{2}(\frac{1}{p}-\frac{1}{q})},
 $$
 where $1\leq p\leq q\leq +\infty$, and $C_p^2= {p^{1/p}} / {p'^{1/p'}}$
 with  $p'$ the  conjugate  index of  $p\in [1,+\infty]$ in the sense that $\frac{1}{p}+\frac{1}{p'}=1$.
This is a particular use  of  the sharp version  of Young's inequality  
for which the  maximizers exist and  are Gaussian functions; see  \cite{LL} p.98.
This is  explicitly written page 223  Equation (2) in \cite{LL}. 
In particular, this yields for  $q=+\infty$, and for all $p$ such that $1\leq p\leq +\infty$,
$$
\vert\vert p_t\star  f
\vert\vert_{\infty}=\vert\vert e^{t\Delta} f
\vert\vert_{\infty}\leq
C(n,p) t^{-\frac{n}{2p}}
\vert\vert f\vert\vert_p, \quad t>0,
$$
where
\begin{equation}\label{cnp}
C(n,p):=C(n,p,+\infty)=
 (p')^{-\frac{n}{2p'}} (4\pi )^{-\frac{n}{2p}}.
\end{equation}

 We use  this contraction property of the heat semigroup on $\R^2$  ($n=2$) as follows. For all $n_0\in L^1(\R^2)\cap L^p(\R^2)$, we have
$$
\vert \vert p_s\star n_0\vert \vert_{\infty}
\leq 
\frac{C(2,p)}{s^{1/p}} \vert \vert n_0\vert \vert_p, \quad s>0,
$$
where $1\leq p\leq +\infty$, and the constant $C(2,p)= (p')^{-\frac{1}{p'}} (4\pi )^{-\frac{1}{p}}$ is obtained from  \eqref{cnp}.
Hence, we deduce that, for all $s>0$ and $z\in \R^2$,
$$
H_z(s)=4\pi s \, p_s\star n_0(z)
\leq 
(4\pi s)^{1-1/p} \, W_p \vert \vert n_0\vert \vert_p,
$$
with $W_p= (p')^{-\frac{1}{p'}}$.
\vskip0.3cm

Now by definition of $T:=T^*_{c,z}(n_0)$, we deduce that
$$
H_{z,n_0}(T)=
L(M)
\leq
(4\pi T)^{1-\frac{1}{p}} W_p \vert \vert n_0\vert \vert_p
=
(4\pi T)^{\frac{1}{p'}} W_p \vert \vert n_0\vert \vert_p.
$$
Since we assume that $p>1$,  we  have  $\frac{1}{p'}=1-\frac{1}{p}>0$.
We easily solve this inequality for $T$,  and get  the following lower bound for $T=T^*_{c,z}(n_0)$, 
$$
\frac{1}{4\pi}
\left[ \frac{L(M)}{W_p \vert \vert n_0\vert\vert_p}
\right]^{p'}
\leq 
T^*_{c,z}(n_0).
$$
By taking the infimum over $z\in \R^2$, we get
$$
\frac{1}{4\pi}
\left[ \frac{L(M)}{ W_p \vert \vert n_0\vert\vert_p}
\right]^{p'}
\leq 
T^*_c(n_0).
$$
Since $n_0\in L^1\cap L^p\subset L^q$ for all $1<q\leq p$, the same  estimate holds with $q$, i.e. 
$$
\frac{q'}{4\pi}
\left[ \frac{L(M)}{  \vert \vert n_0\vert\vert_q}
\right]^{q'}
=
\frac{1}{4\pi}
\left[ \frac{L(M)}{ W_q \vert \vert n_0\vert\vert_q}
\right]^{q'}
\leq 
T^*_c(n_0).
$$
Finally, we conclude the proof of \eqref{lowerbd} by taking the supremum over $q\in (1, p]$.
\vskip0.3cm

(2)\,  Let $1<p\leq +\infty$.
By the generalized H\"older's inequality, we can write for all $1<q\leq p$,
$$
 \vert \vert n_0\vert\vert_q
 \leq
 M^{1-p'/q'}  \vert \vert n_0\vert\vert_p^{p'/q'},
 $$
where $M= \vert \vert n_0\vert\vert_1$.
Here, $p^{\prime},q^{\prime}$ are respectively the  indices conjugate to $p,q$.
From  the inequality  \eqref{lowerbd},
we deduce that
$$ 
\frac{1}{4\pi}
\sup_{1<q\leq p} q'
\left[ \frac{L(M)}{  M^{1-p'/q'}  \vert \vert n_0\vert\vert_p^{p'/q'} }
\right]^{q'}
\leq 
T^*_c(n_0),
$$ 
or equivalently,
$$
\frac{1}{4\pi}
\sup_{q'\in [p', \,+\infty)}
 q'
\left[ \frac{L(M)}{M}\right]^{q'}
 \left[\frac{M} {\vert \vert n_0\vert\vert_p}\right]^{p'} 
\leq 
T^*_c(n_0).
$$
This lower bound on $T^*_c(n_0)$ can be written as,
$$ 
\frac{1}{4\pi} \left[ \frac{M} {\vert \vert n_0\vert\vert_p}\right]^{p'} 
\sup_{q'\in [p',\,+\infty)}
 q'
a^{q'}
\leq 
T^*_c(n_0).
$$ 
where $a:=\frac{L(M)}{M}=\frac{2M}{3M-8\pi}$. 
\vskip0.3cm

Thus, we just need to estimate the supremum just  above.
Set $h(r)=r a^r$, $r>0$. Recall that $a\in (0,1)$ when $M>8 \pi$.
It is easily shown that $h$ is increasing on the interval $(0,r_0)$
and decreasing on $(r_0,+\infty)$,  where 
$$
r_0=\frac{1}{\ln (1/a)}=
\frac{1}
{\ln\left(1+\frac{M-8\pi}{2M}\right)}.
$$
Thus,  $h$ attains its maximum at $r_0$ on $(0,+\infty)$. We deduce  that
$$
\sup_{r>0} h(r)={\rm max}_{r>0} \,  h(r)=
h(r_0)=
\frac{e^{-1}}
{\ln\left(1+\frac{M-8\pi}{2M}\right)}.
$$
This leads us to discuss two cases for evaluating  $\sup_{q'\in [p'\,+\infty)} q'a^{q'}$ for    fixed $p'\in [1,+\infty)$.
\vskip0.3cm 

(i)
If $ p' \leq r_0$, then we have
$$
\sup_{q'\in [p'\,+\infty)}
 q'
a^{q'}
=\frac{e^{-1}}
{\ln\left(1+\frac{M-8\pi}{2M}\right)}.
$$
\vskip0.3cm

(ii)\, 
If $p' >r_0$, then we have
$$
\sup_{q'\in [p'\,+\infty)}
 q'
a^{q'}
=
 p'
a^{p'}.
$$
\vskip0.3cm 

This  implies that
\vskip0.3cm

(iii)\,
If $3-2e^{1-1/p}\leq  \frac{8\pi}{M}$
(i.e. $ p' \leq r_0$),
then we have
$$ 
\frac{(\pi e)^{-1}}
{4 \ln\left(1+\frac{M-8\pi}{2M}\right)}
\left[ \frac{M}{  \vert \vert n_0\vert\vert_p}
\right]^{p'}
\leq 
T^*_c(n_0).
$$

(iv)\, 
If  $3-2e^{1-1/p}> \frac{8\pi}{M}$
(i.e. $ p'>r_0$),
 then we have
$$ 
\frac{p'}{4\pi} \left(\frac{L(M)} {\vert \vert n_0\vert\vert_p}\right)^{p'} 
\leq 
T^*_c(n_0).
$$
\vskip0.3cm

Now we are in position to discuss both cases (a) and (b)  of Theorem \ref{lowe1} (2).
\vskip0.3cm

(a) Let $p_0=\left[ \ln (\frac{2e}{3})\right]^{-1}\; ( \sim 1,682)$. Assume that $p_0\leq p\leq +\infty$. 
This condition is equivalent to $3-2e^{1-1/p}\leq 0$ ($p_0$ is the unique  root of $v(p)=3-2e^{1-1/p}$).
 Thus, the condition  of  (iii) just above
 is trivial  for any $M>8\pi$. Consequently, the inequality \eqref{lowerbd2} holds true.
 Now, let $1<p<p_0$ and $M\leq \frac{8\pi}{3 -2e^{1-1/p}}$. 
 Then this last  inequality  is equivalent to
$0<3-2e^{1-1/p}\leq \frac{8\pi}{M}$, and again \eqref{lowerbd2} holds true by virtue of  (iii).
\vskip0.3cm

(b) Let $1<p<p_0$ and $M >  \frac{8\pi}{3 -2e^{1-1/p}}$.  Then this second inequality  is equivalent to
$3-2e^{1-1/p} > \frac{8\pi}{M}$, hence the inequality \eqref{lowerbd3} holds true by virtue of (iv).

\vskip0.3cm

This concludes the proof of  Theorem  \ref{lowe1}.
\hfill $\square$
\vskip0.3cm

\begin{rem}
Let $a=\frac{L(M)}{M}=\frac{2M}{3M-8\pi}$ with $M>8\pi$. From the study of the function  $h$ defined  above,  
it is easy to compare  formally both inequalities  \eqref{lowerbd2} and \eqref{lowerbd3}. 
We can see that  \eqref{lowerbd2}  is   stronger than  \eqref{lowerbd3} (at least for $p'\leq r_0$).
We have the next inequality,
$$
h(p')=p' a^{p'} 
<
\sup_{r\in [p',+\infty)}
h(r)=h(r_0)=
\frac{ e^{-1}}
{ \ln\left(\frac{1}{a} \right)}.
$$
Now, assuming \eqref{lowerbd2}
we  would have the next inequality,
$$ 
\frac{p'}{4\pi} \left(\frac{L(M)} {\vert \vert n_0\vert\vert_p}\right)^{p'} 
=
\frac{p'}{4\pi}{a^{p'}} \left(\frac{M} {\vert \vert n_0\vert\vert_p}\right)^{p'} 
\frac{ e^{-1}}
{4 \pi \ln\left(\frac{1}{a} \right)}
\left(\frac{M}{  \vert \vert n_0\vert\vert_p}
\right)^{p'}
\leq
T^*_c(n_0),
$$
i.e.  \eqref{lowerbd3}. This proves the assertion.
\end{rem}

\section{Examples of initial data}\label{examples}
This section is devoted to applications of    Theorem \ref{tcalphazero}, \ref{tclaplace},  and 
 Corollaries  \ref{corapplica0}, \ref{corapplica1}, \ref{corapplica3}, \ref{corapplica2}, providing upper  estimates on the critical time  $T^*_c(n_0)$
to solutions of the  (PKS)   system for  several families  of explicit  initial data $n_0$ with supercritical mass $M>8\pi$.
Theorem \ref{lowe1} is  applied to provide lower bounds on  $T^*_c(n_0)$.
We shall also compare the different  estimates of  $T^*_c(n_0)$ obtained  from these methods.
As a consequence, we deduce explicit upper bounds for the  maximal existence time $T^*$ of the (PKS) system.
 See Sections \ref{intro}, \ref{tcsection} and \ref{lowerboundsect}.
Note that  Corollary \ref{corapplica0} is not used  in this section for  the examples described  below due to the difficulty 
for  obtaining explicit computations.

\subsection{Gaussian initial data $n_0$}
Our first family of examples of initial data $n_0$ consists of  Gaussian functions written in the following form,
$$
n_0(x)=
\frac{M}{4\pi\sigma}e^{-\frac{\vert x-z_0\vert^2}{4\sigma}}
=Mp_{\sigma}( x-z_0),\quad x\in \R^2,
$$
where $\sigma>0$ and $z_0\in \R^2$ are fixed. 
Here, the function $p_{\sigma}$ denotes the  Gaussian (or  heat)  kernel.
For all these data $n_0$,
the mass is $\vert\vert n_0\vert\vert_1=M$
and the normalized barycenter is  $B_0=z_0$.
In the whole section, we  shall assume that $M>8\pi$. 
By applying part (3) and (4) of Theorem \ref{tcalphazero}, 
we obtain the exact value of $T^*_c(n_0)$,
and  also the next upper bound for the maximal existence time $T^*$ of the (PKS) system,
\begin{equation}\label{gaussianexam}
T^*
\leq
T^*_c(n_0):=\sigma\frac{2M}{M-8\pi}.
\end{equation}
The proof is as follows. The function  $n_0$ is a non-increasing
$z_0$-radially symmetric function. We set $m_0(x)=n_0(x+z_0)$, $x\in \R^2$. 
So, $m_0$  is a non-increasing radially symmetric function.
By (4) of Theorem \ref{tcalphazero}, we have 
$T^*_c(n_0)= T^*_c(m_0)$.
We determine  the function $H_{z,m_0}$ evaluated at $z=0$ and get  the next formula 
$$
H_{0,m_0}(s)=4\pi sM \, p_s\star p_{\sigma}(0)
=4\pi sM \, p_{s+{\sigma}}(0)
=
\frac{sM}{s+\sigma}, \; s>0.
$$
The inverse of $H_{0,m_0}$ is easily deduced,   we have  
$$
s= H^{-1}_{0,m_0}(u)=
\frac{\sigma u}{M-u}, \; u\in (0,M)
.
$$
Finally,  we deduce  by (3) of  Theorem \ref{tcalphazero} that
$$
T^*_c(n_0)=T_c^*(m_0)=H^{-1}_{0,m_0}\left( \frac{2M^2}{3M-8\pi}\right)
=\sigma\frac{2M}{M-8\pi}.
$$
Alternatively, we can apply  \eqref{laplacegal1} of Theorem \ref{tclaplace} (but the computation is longer) 
and obtain the same result for the value of $T^*_c(n_0)$.
Note that $T^*_c(n_0)$ depends not only of the mass $M$ but also of the {\em shape} of the Gaussian $p_{\sigma}$ via the function
$H^{-1}_{0,m_0}$, in particular
through the variance parameter $\sigma$.
Of course, the {\em shape} does not depend on the spacial  {\em position} parametrized by  the translation term $z_0$ 
of the Gaussian function $p_{\sigma}$. In the limit case of concentration $\sigma=0$, we would  obtain
$T^*_c(n_0)=0$ (the mass $M$ remains fixed), i.e. an {\em immediate blow-up} of the solution of the (PKS) system as it can be expected.
This  is not physically surprising since the initial data should be  considered at this limit case as a Dirac measure at $z_0$
and the result  is consistent with our intuition. On the other hand, for the behavior of a radially symmetric  solution (i.e. $z_0=0$)
 closed to the blow-up time $T^*$, we can refer for instance  to \cite{BDP, BCM} (and the references therein) for a mathematical statement of the blow-up profile.
\vskip0.3cm

As already mentioned, whenever the second moment of $n_0$ is finite and for which an estimate is available, 
we shall compare  the upper bound  $T^*_v$ of  $T^*$ (see \eqref{bupbdp2}) and  the one given by $T^*_c(n_0)$ obtained in this paper.
This is  done in (1) below.
\vskip0.3cm

\begin{enumerate}
\item
 {\em  Second moment method}.
We deduce from  the inequality  \eqref{bupbdp2}
 the next estimate of $T^*$,
$$
T^*\leq T^*_v:=\frac{2\pi V_2(n_0)}{M-8\pi}=
\sigma\frac{8\pi}{M-8\pi},
$$
for all $\sigma>0$ and all $z_0\in \R^2$.
Indeed, it is easy to check that $B_0=z_0$ and the well-known variance of the Gaussian function $p_{\sigma}$ in $\R^2$
is given by
\begin{equation}\label{var2}
V_2(n_0)=
\frac{1}{M}\int_{\Ri^2} \vert x- B_0\vert^2 n_0(x)\, dx
=\int_{\Ri^2} \vert x\vert^2 p_{\sigma}(x)\, dx= 4\sigma .
\end{equation}
This upper bound $T^*_v$  is obviously a better bound  of $T^*$ than the one given by \eqref{gaussianexam} since $M>8\pi$. 
This is due to the fact that $T^*_v$ is deduced directly  from the second moment evolution equation for a solution of the (PKS) system
 under finite second moment   assumption of the initial data $n_0$ (see details below   \eqref{bupbdp2}).
This  additional information on $n_0$ leads to  a better estimate  compared to  the general bound  $T^*_c(n_0)$ in  \eqref{tcnormal}   of $T^*$ 
valid for any initial data $n_0$.
\vskip0.3cm

\item
 {\em Various estimates from the general case}.
As much as possible, we now  explicit  various estimates of $T^*_c(n_0)=\sigma\frac{2M}{M-8\pi}$  obtained by the  corollaries  described after  
Theorem \ref{tcalphazero} and by Theorem \ref{tclaplace}, and compare them with the exact  formula \eqref{gaussianexam} of $T^*_c(n_0)$.
 Note that this comparison concerns the upper bound $T^*_c(n_0)$ of $T^*$ and not $T^*$ itself.
 Below, we present the results not necessarily in the order  of the corollaries.
  \end{enumerate}
\vskip0.3cm

{\em About lower bounds given by  Theorem \ref{lowe1}}.  The inequality   \eqref{lowerbd}  
is, in fact, an equality, i.e.
\begin{equation}\label{loweroptim}
\frac{1}{4\pi}
\sup_{1<q\leq p} q'
\left[ \frac{L(M)}{  \vert \vert n_0\vert\vert_q}
\right]^{q'}
=
\sigma\frac{2M}{M-8\pi}
=
T^*_c(n_0).
\end{equation}
It comes as no surprise  because the proof of  Corollary \ref{lowe1} relies on   sharp Young's inequalities for which optimizers are Gaussians
(see \cite{LL} p.99).
This shows  that  the inequality \eqref{lowerbd}  can be  optimal at least for  this family of examples.
Let us prove this equality.
For  $1< q <+\infty$, the $L^q$-norm of $n_0$ is given explicitly by
$$
\vert \vert n_0\vert\vert_q=
M \vert \vert p_{\sigma}\vert\vert_q
=M\frac{(4\pi\sigma)^{-1+1/q}}
{q^{1/q}}.
$$
Thus, we can set 
$$
{\mathcal R}(q'):=
\frac{ q'}{4\pi}
\left[ \frac{L(M)}{  \vert \vert n_0\vert\vert_q}
\right]^{q'}
=\sigma a^{q'} q'q^{q'/q},
$$
with $a=\frac{L(M)}{M}=\frac{2M}{3M-8\pi}$. Now, since 
$q'/q=q'-1$ and $q=\frac{q'}{q'-1}$, we can write 
$$
{\mathcal R}(q')=\sigma a^{q'} q'q^{q'-1}
=
\sigma a^{q'} (q'/q) q^{q'}
=
\sigma a^{q'} (q'-1) \left( \frac{q'}{q'-1}\right)^{q'}=
\sigma a^{q'} (q'-1)^{1-q'} (q')^{q'}.
$$
We set $\rho=q'\geq p'$. Then, for $p=p(M)$ large enough so that $p'\leq \frac{1}{1-a}$, we have
$$
\frac{1}{4\pi}
\sup_{1<q\leq p} q'
\left[ \frac{L(M)}{  \vert \vert n_0\vert\vert_q}
\right]^{q'}
=
\sup_{q'\geq p'} 
{\mathcal R}(q')
=
\sigma  \sup_{q'\geq p'} 
a^{q'} (q'-1)^{1-q'} (q')^{q'}
$$
$$
=
\sigma  \sup_{\rho\geq p'} 
e^{H(\rho)}
=\sigma  \sup_{\rho\geq 1} 
e^{H(\rho)}
=
\sigma e^{H(\frac{1}{1-a})}
=
\sigma \frac{a}{1-a}
= \sigma\frac{2M}{M-8\pi},
$$
where 
$H(\rho)=  \rho \ln a+(1-\rho)\ln(\rho -1)+\rho \ln \rho $.
Indeed, from the study of the derivative 
$$
H'(\rho)=
\ln a +\ln \left(\frac{\rho}{\rho-1}\right),
$$
we deduce that the function $\rho\mapsto H(\rho)$ is  increasing on the interval  $(1,\frac{1}{1-a})$, and decreasing on  
 the interval  $(\frac{1}{1-a},+\infty)$.
 Thus, the supremum  of $H$ on $(1,+\infty)$ is achieved for $\rho_0=\frac{1}{1-a}>1$. 
Now, since $n_0$ is in $L^p$ for all $p\geq1$, 
we can choose $p$ large enough,  so that $p'\leq \rho_0$.
 This leads to the expected equality \eqref{loweroptim}.
\vskip0.3cm

The estimates  of $T^*_c(n_0)$ obtained from part 2 of  Theorem \ref{lowe1} are, of course,  less precise.
Just below, we list them. Under the corresponding assumptions for each assertion, we can write
\vskip0.3cm

\noindent
(a)
 from \eqref{lowerbd2}, 
$$
 \sigma \frac{e^{-1} p^{p' -1}}{
\ln (1+\frac{M-8\pi}{2M})}\leq   T^*_c(n_0),
$$
(b)
from \eqref{lowerbd2bis}, 
$$
\sigma \frac{e^{-1} }{
\ln (1+\frac{M-8\pi}{2M})}\leq   T^*_c(n_0),
$$
(c)
 from  \eqref{lowerbd3},  
 $$
\sigma \frac{p'\, p^{p'-1} }{
(1+\frac{M-8\pi}{2M})^{p'}}\leq   T^*_c(n_0).
$$
A priori, none of  these inequalities are  sharp.
They don't lead to the optimal estimate \eqref{lowerbd}  of $T^*_c(n_0)$ (case of equality for Gaussians). 
Nevertheless,  the inequalities \eqref{lowerbd2}  and \eqref{lowerbd2bis}  are all of the  next form,
$$
 \sigma \,  \frac{k}{
\ln (1+\frac{M-8\pi}{2M})}\leq   T^*_c(n_0),
$$
for  some  constant $k>0$ (independent of $M$). 
Then we can   compare these estimates  with the exact value of  $T^*_c(n_0)$  as $M$ goes to $8\pi^+$,
$$
 \sigma \frac{ k}{
\ln (1+\frac{M-8\pi}{2M})}\ \sim  k 
 \frac{2 {\sigma}M}{M-8\pi}=k T^*_c(n_0).
$$
\vskip0.3cm

We now  compare the  exact value of  $T^*_c(n_0)$  with various upper bounds on $T^*_c(n_0)$  
obtained  from  the corollaries of Theorem \ref{tcalphazero} of Section \ref{tcsection} 
for  Gaussian initial data.
\vskip0.3cm

{\em  From Corollary  \ref{corapplica1}}. 
In order to  apply  \eqref{radial},  
we  first  compute
$$
{\mathcal M}_{z_0}(\rho):=  
\int_{B(z_0,\rho)} n_0(x)\, dx=
M(1-e^{-\frac{\rho^2}{4\sigma}}), \; \rho>0.
$$ 
Hence, we get 
$$
\frac{
{\mathcal M}_{z_0}(\rho)
}
{L(M)}
= \frac{1}
{a}( 1-e^{-\frac{\rho^2}{4\sigma}}),
$$
where $\frac{1}{a}=1+\frac{M-8\pi}{2M}$. Using the fact that $n_0$ is  non-increasing  
$z_0$-radially symmetric, we can apply    \eqref{radial} (by invariance, it is sufficient  to  consider the case $z_0=0$).
Thus we obtain 
$$
T_c^*(n_0) \leq T_{c_2}^*
:=
\inf_{\rho>0}
\frac{\rho^2}{4}\left[\ln_+\left(\ \frac{1-e^{-\frac{\rho^2}{4\sigma}}}
{a}
\right)\right]^{-1}.
$$
By setting  $s=\frac{\rho^2}{4\sigma}$,  this implies that
$$
T_c^*(n_0) \leq T_{c_2}^*:=
\sigma
\inf_{s>0}
{s}
 \left[\ln_+\left(\ \frac{1-e^{-s}}
{a}
\right)\right]^{-1}
=
\sigma
\inf_{s>0,1-e^{-s}>a}\,
\frac{s}
{
\ln (1-e^{-s})+\ln(1/a)}.
$$
Unfortunately, it is not clear if we can estimate this infimum in such  a way that it is closed to the value of $T_c^*(n_0)$ given above.
We can   deduce  that we have for all $\varepsilon \in (0,1-a)$, 
$$
T_c^*(n_0) 
 \leq \sigma
 \frac{ c_{\varepsilon}}
{\ln(\frac{1-\varepsilon}{a})},
$$
with  $c_{\varepsilon}=\ln(\frac{1}{\varepsilon})>0$ by taking $s=c_{\varepsilon}$.
 We shall see below that we can obtain a more interesting   upper bound  on $T_c^*(n_0) $ of the form 
 $\sigma \frac{1} {\ln(\frac{1}{a})}$ using the finite variance of the initial data $n_0$
  by applying Corollary  \ref{corapplica2}. This will improve the upper bound on $T_c^*(n_0) $ just above obtained by  Corollary  \ref{corapplica1}.
\vskip0.3cm

Another bound is also obtained from the  expression of  ${\mathcal M}_{z_0}(\rho)$. 
Indeed, we deduce that
$$
g_{z_0}(\rho):=\frac{1}{M} \int_{B(z_0,\rho)} n_0(x)\, dx= 
1-e^{-\frac{\rho^2}{4\sigma}}.
$$
Thus, for all  $ m\in (0, 1)$ we have
$$
h(m)
=\left[
g_{z_0}^{\leftarrow}(m)
\right]^2
=\left[
g_{z_0}^{-1}(m)
\right]^2
=\rho^2=-4\sigma \ln (1-m).
$$
Applying \eqref{tdeuxd}, we can write
$$
T^*_{c_2}
\leq
\frac{\sigma}
{
\ln\left(
1+\frac{M-8\pi}{2M}
\right)
}.
\inf_{\theta \in (0,1)}
\left(
\frac{ 
-\ln(1-a^{\theta})
}{(1-\theta)}
\right).
$$
Here again, it is not clear   if we can estimate this infimum in such a way that the right-hand side of the preceding inequality  is closed to $T_c^*(n_0)$.
\vskip0.3cm

Finally, note that the inequalities \eqref{tctrois} and \eqref{tctroisbis} of  Corollary  \ref{corapplica1} do not apply 
 here since $n_0$ has no compact support.
\vskip0.3cm

{\em  From   Corollary  \ref{corapplica3}}.
We can compute explicitly  $F$ and study its properties. We have
$$
\frac{h}{h'}(m)=-(1-m)\ln(1-m), \; m\in (0,1).
$$
Then $F$ is given by 
$$
F(X)=X+e^X \left(\frac{h}{h^{\prime}}\right)(e^{-X})
=
X-(e^X-1)\ln (1-e^{-X}), \quad X>0.
$$ 
and its derivative by
$$
F^{\,\prime}(X)=
-e^X \ln (1-e^{-X}), \quad X>0.
$$
The condition \eqref{tcF1}, namely  $ \left(\frac{h}{h^{\prime}}\right)(1^-)=0\geq \ln(1/a)$, is not satisfied since $a\in (0,1)$. 
Then (1) of  Corollary  \ref{corapplica3} does not apply.
But the conditions of part (2) are  easily satisfied. Indeed, 
the function $F$ is clearly continuous and strictly increasing  ($F^{\,\prime}(X)>0$ for all $X>0$). 
Condition (i) is satisfied, i.e. 
 $ \left(\frac{h}{h^{\prime}}\right)(1^-)=0< \ln(1/a)$ holds true, and 
 (ii) is also  satisfied, i.e. 
 $$
  \left(\frac{h}{h^{\prime}}\right)(a^+)=-(1-a) \ln(1-a)>0
 .
 $$
By applying  \eqref{tcF2}, we get
$$
T^*_{c_2}(0)=
S\left(Y_0\right)
=\frac{\sigma}{e^{F^{-1}(Y_0)}-1},
$$
where $Y_0=\ln \frac{1}{a}$. 

We now compare $T^*_{c_2}(0)$ with  $T^*_{c}(n_0)$.
It is clear that for all $X>0$, we have $X<F(X)$. So, we get 
$F^{-1}(Y) <Y$  for all $Y\in {\rm Im} (F)$. This implies that
$$
\frac{\sigma}{e^Y-1},
<
\frac{\sigma}{e^{F^{-1}(Y)}-1}
$$
which, evaluated at $Y=Y_0$, gives us
$$
 T^*_{c}(n_0)=  
\sigma\frac{2M}{M-8\pi}
=
\frac{\sigma}{e^{Y_0}-1}
<
\frac{\sigma}{e^{F^{-1}(Y_0)}-1}=
T^*_{c_2}(0).
$$
This proves again that $ T^*_{c}(n_0)\leq T^*_{c_2}(0)$ but  clearly 
 $T^*_{c}(n_0)\neq T^*_{c_2}(0)$ for the family  of Gaussian  initial data.
 
The inverse function $F^{-1}$ of $F$ can not be given explicitly for  this particular case.
Then  we look for  an upper bound on $ T^*_{c_2}(0)$.
A first approach should be  to  note  that $F(X)\leq X+G(\ln(1/a))<X+1$ for $X\in (0,\ln(1/a))$ with  $G(X)=F(X)-X$ ($G$ is an increasing function).
This  yields  $Y_0-G(\ln(1/a))\leq F^{-1}(Y_0)$. But unfortunately $Y_0-G(\ln(1/a)<0$.
Thus we cannot deduce an upper bound for $T^*_{c_2}(0)$ from this approach.
Another approach will be  to find  some finite  constant $c>0$ such that  $F(x)\leq cX$  for all $X\in (0,\ln(1/a))\subset (0, 0.41)$.
But we are facing to the problem that  $\lim_{X\rightarrow 0^+} F(X)/X=+\infty$  which contradicts such an estimate.
Since we cannot conclude  here on  this point,  we shall not go further in our analysis of $F^{-1}$.
\vskip0.3cm

{\em From  Corollary  \ref{corapplica2}}.
We obtain an  explicit upper bound on $T^*_c(n_0)$ using  $ T^*_{c_4}$ given by \eqref{tstarparti2}.
This estimate  is sharp when $M$ is closed to $8\pi$.
We also deduce  interesting upper and lower bounds  of $T^*_c(n_0)$ in terms of  $ T^*_{c_4}$.
First, from the value of the variance  $V_2(n_0)$ of $n_0$ given by \eqref{var2} and  applying \eqref{tstarparti2}, it follows  immediately that
$$
T^*_c(n_0)=\sigma\frac{2M}{M-8\pi}
\leq T^*_{c_4}=\sigma\frac{1}{
\ln (1+\frac{M-8\pi}{2M})
}.
$$
This explicit  result  is  sharp  when $M\rightarrow 8\pi^+$, since we have the following  asymptotic behaviour
$$
 T^*_{c_4}=\sigma\frac{1}{
\ln (1+\frac{M-8\pi}{2M})}  \sim
\sigma\frac{2M}{M-8\pi}=T^*_c(n_0),  \quad  M\rightarrow 8\pi^+.
$$
On the other hand, by a  direct study of the  function  $y\mapsto v(y)=\frac{\ln(1+y)}{y}$ which  is  a decreasing function on $(0,1/2)$,
we deduce that
$$
c_0\frac{1}{\ln(1+y)} \leq \frac{1}{y} \leq \frac{1}{\ln(1+y)}, \quad y\in (0,1/2),
$$
where $c_0:=v(1/2)={2\ln (3/2)}$.
We apply this inequality to  $y=\frac{M-8\pi}{2M}\in (0,1/2)$ with $M>8\pi$, and get 
$$
{c_0} \frac{\sigma}{
\ln (1+\frac{M-8\pi}{2M})}
\leq
\sigma \frac{2M}{M-8\pi}
\leq
\frac{\sigma}{
\ln (1+\frac{M-8\pi}{2M})}, \quad M>8\pi.
$$
In other words, we have obtained  the following comparison  estimates between $T^*_c(n_0)$ and  $T^*_{c_4}$ ,
$$
 {c_0}  T^*_{c_4} \leq  T^*_c(n_0)\leq T^*_{c_4},
$$
where  ${c_0}  \sim 0. 8109$ and  ${c_0}\geq  0. 8109$.
Note that the upper  bound  just above is already known from Corollary  \ref{corapplica2}.
The lower bound is of interest for estimating $T^*_c(n_0)$ by keeping  in mind the fact that
$T^*_c(n_0)$ tends to infinity  when  $M\rightarrow 8\pi^+$  (hence also $ T^*_{c_4}$). 
\vskip0.3cm

Recall here   that our main goal in this section  is to compare general results obtained in Corollaries
\ref{corapplica1}, \ref{corapplica3}, \ref{corapplica2} and Theorem \ref{lowe1} with the known and explicit value
$$
  T^*_{c}(n_0)
  =
  \sigma \frac{2M}{M-8\pi},
  $$
obtained by Theorem \ref{tcalphazero}, or Theorem \ref{tclaplace}.
We  summarize   the most interesting results of this comparison of $T^*_c(n_0)$ with various $T^*_{c_i}$
for the gaussian  initial data as follows.
 From Corollary  \ref{corapplica2} using the variance of the initial data and Theorem \ref{lowe1}, 
 we have obtained  the following  estimates of $T^*_{c}(n_0)$, namely
$$
 T^*_{c_5}=
  T^*_{c}(n_0)= \sigma \frac{2M}{M-8\pi}
  \leq
  T^*_{c_4}=
   \sigma \frac{1}{
\ln (1+\frac{M-8\pi}{2M})}
\leq
 \frac{1}{c_0} T^*_{c}(n_0).
$$
Moreover,  Theorem \ref{lowe1} provides  an optimal result (equality case).
We have  also seen that the estimate of $T^*_{c}(n_0)$ by $T^*_{c_4}$ of  Corollary  \ref{corapplica2}
 is asymptotically sharp in the sense that
$T^*_{c_4}  \sim  T^*_{c}(n_0) $ when $M\rightarrow 8\pi^+$. 
The method of estimating $T^*_{c}(n_0)$ from the variance is particularly interesting for the Gaussian initial data.

\subsection {Characteristic function  of a disk} 
A natural  situation  is when the cells are uniformly distributed in a disk $B(z_0,R)$ with some accumulation height  $\sigma>0$. 
This leads to  consider an initial data of the type  $n_0=\sigma {\bf 1}_{B(z_0,R)}$ for fixed $z_0\in \R^2$ and $R>0$.
Here $ {\bf 1}_{B(z_0,R)}$ denotes the characteristic function  of the disk ${B(z_0,R)}$ of radius $R$ centered at  $z_0$.
 We assume to be in  the supercritical case  $M>8\pi$,
i.e.  $ \sigma R^2>8$. In that case, we can describe  the bound $T^*_c(n_0)$ as follows.

\begin{pro}\label{ball}
Let $n_0=\sigma {\bf 1}_{B(z_0,R)}$ be the initial data for the (PKS) system.
Here,  $z_0\in \R^2$ and $\sigma>0$ are fixed.
Assume that $M=\pi \sigma R^2>8\pi$. 
\begin{enumerate}
\item
Then the maximal  existence  time $T^*$ of the  solution of  (PKS) system with initial data $n_0$
is bounded as follows:
\begin{equation}\label{tcbound1}
T^*\leq T^ *_c(n_0):=\frac{R^2}{4f^{-1}(\frac{2M}{3M-8\pi})}
=
\frac{M}{4\pi \sigma f^{-1}(\frac{2M}{3M-8\pi})}
,
\end{equation}
where $f^{-1}$ is the inverse of the  function
$f(\lambda)=\frac{1-e^{-\lambda}}{\lambda}$, $\lambda>0$.
In particular, we have
\begin{equation}\label{tcbound2}
0,286.R^2\sim
\frac{R^2}{4f^{-1}(\frac{2}{3})}\leq T^ *_c(n_0),
\end{equation}
with $f^{-1}(\frac{2}{3})\sim 0.874 \,21$.
As a consequence of \eqref{tcbound1}, let  $\sigma \rightarrow +\infty$ for fixed $R>0$  (so, $M \rightarrow +\infty$), we have
\begin{equation}\label{tcbound3}
T^ *_c(n_0)\sim \frac{R^2}{4f^{-1}(\frac{2}{3})}\sim0,286.R^2.
\end{equation}
\item
We have also the following uniform estimate for $M>8\pi$,
\begin{equation}\label{tcbound2bis}
e^{-1} 
\frac{R^2}
{4\ln\left(1+\frac{M-8\pi}{2M}\right)}
\leq T^ *_c(n_0)
\leq
\frac{R^2}{4\ln\left(1+\frac{M-8\pi}{2M}\right)}.
\end{equation}
\item
The following asymptotic estimate  holds true   when $M$ is closed to the supercritical mass $8\pi$. 
\begin{enumerate}
\item
Let $R>0$ be fixed. Then we have
\begin{equation}\label{tcbound4}
T^ *_c(n_0)\sim \frac{2\pi R^2}{M-8\pi}=\frac{2 R^2}{\sigma R^2-8}, \; \; {\rm as} \;  \;  \sigma \rightarrow\left(8/R^2\right)^+.
\end{equation}
\item
Let $\sigma>0$ be fixed. Then we have
\begin{equation}\label{tcbound5}
T^ *_c(n_0)\sim \frac{16\pi }{\sigma(M-8\pi)}=\frac{16}{\sigma(\sigma R^2-8)}, \; \; {\rm as} \; \;  R \rightarrow\left(\sqrt{8/\sigma}\right)^+.
\end{equation}
\end{enumerate}
\end{enumerate}
\end{pro}

As already mentioned,  the bound  $T^ *_c(n_0)$ does not depend  on  the normalized barycenter $B_0=z_0$ by translation invariance.
Indeed, we should expect from a physical point of view that the evolution of cells  should be the same if we make a translation
 in the space variable  by any $z_0\in \R^2$ of the initial data $n_0$ by considering
the new initial data  $x\mapsto n_0(x+z_0)$. Of course, this  due to the implicit isotropic environment of $\R^2$ here (flat curvature).
\vskip0.3cm

As described by   \eqref{tcbound1},  the critical time $T^ *_c(n_0)$ can  be  alternatively expressed   with  the radius $R$ 
(which measures essentially the size of the support of $n_0$),
or   with the height $\sigma$ (i.e. the sup norm  of $n_0$) for  $n_0$ with fixed mass  $M>8\pi$.
 As a consequence, if we fix the mass $M>8\pi$ and  let $R\rightarrow 0$ (i.e., $\sigma \rightarrow +\infty$), we deduce that $T^ *_c(n_0)$ tends to $0$,
and consequently  $T^ *$ too.
In particular, again  with a  fixed mass $M>8\pi $, if we consider a sequence of initial data  of the form $n_{0,k}=\sigma_k {\bf 1}_{B(z_0,R_k)}$  
with $R_k \rightarrow 0$ 
(i.e., $\sigma_k=\frac{M}{\pi R_k^2}\rightarrow +\infty$), which approximates $M\delta_{z_0}$ in a weak sense,
then we obtain from \eqref{tcbound1} that $\lim_k T^*_c(n_{0,k})=0$.   
Here,  $\delta_{z_0}$ denotes the Dirac measure  at $z_0\in \R^2$.
Heuristically, it should   not be surprising to observe  an  instantaneous  explosion if we consider   $M\delta_{z_0}$ as initial data with $M>8\pi$.
\vskip0.3cm

Now if we fix  $R>0$ (i.e. the support of $n_0$ is fixed)  and let the accumulation height 
$\sigma$ of cells goes to infinity, then we deduce from \eqref{tcbound3}
 that  the critical time $T^ *_c(n_0)$, respectively $T^*$, is uniformly bounded above  by $cR^2$ for some  constant $c$.
\vskip0.3cm
 
{\bf Proof of Proposition \ref{ball} }.
(1)\,
We start from the upper bound  \eqref {heatcondit} on $T^*$. Since the function 
$n_0(x)=\sigma {\bf 1}_{B(z_0,R)}(x)=\sigma {\bf 1}_{B(0,R)}(x-z_0)$  
is the translated by $z_0$ of a radial function, then we have 
$$
T^*_c(n_0)=H_{z_0}^{-1}\left(\frac{2M^2}{3M-8\pi}\right),
$$
where
$$ 
 H_{z_0}(T)=
 \int_{\Ri^2} \exp\left(-\frac{\vert x -z_0\vert^2}{4T}\right)n_0(x)\, dx
 =
\sigma  \int_{B(0,R)} \exp\left(-\frac{\vert y\vert^2}{4T}\right)\, d 
 y
$$
 $$
 =
 2\pi \sigma 
 \int_0^R
 \exp\left(\frac{-r^2}{4T}\right)\, rdr 
  =
 2\pi\sigma 
 \left[-2T\exp\left(\frac{-r^2}{4T}\right)\right]_{r=0}^{r=R}
 =
  4\pi\sigma T(1- e^{-\frac{R^2}{4T}}).
 $$
This result  is obtained by the change of variables $y=x-z_0$ and using polar coordinates.
Now from the relation $\pi \sigma=\frac{M}{R^2}$,  this implies
$$
H_{z_0}(T)=  M \left(\frac{4T}{R^2}\right)
\left(1-e^{-\frac{R^2}{4 T}}
\right)
=M f \left(\frac{R^2}{4T} \right),
$$
where the function $f:(0,+\infty)\rightarrow (0,1)$ is defined by
$f(\lambda)=\frac{1-e^{-\lambda}}{\lambda}$.
\vskip0.3cm

Since the function $f$  is a continuous (convex)  strictly decreasing function with range $(0,1)$ then $f$ is invertible. 
Then the  inverse function $f^{-1}:(0,1) \rightarrow (0,+\infty)$  of $f$  is  also a continuous (convex)  strictly decreasing function.
Hence, we  obtain  the following expression
$$
H_{z_0}^{-1}(u)=
\frac{R^2}{4f^{-1}(\frac{u}{M})},
$$
for all  $u\in (0,M)$.
From the relation $M=\sigma \pi R^2$, this  yields the next upper bound on $T^*$,
$$
 T^*
 \leq
 T^*_c(n_0):=
 H_{z_0}^{-1}
 \left( 
 \frac{2M^2}{3M-8\pi}
 \right)
 =
\frac{R^2}{4f^{-1}(\frac{2M}{3M-8\pi})}
=
 \frac{M}
 { {4\pi\sigma}
 f^{-1}\left(
\frac{2M}{3M-8\pi}
\right)}.
$$
This proves  the inequality \eqref{tcbound1}  which is the first part of (1) of  the proposition.
\vskip0.3cm

We  finish the proof of this  part (1) as follows.
First, we note that  $\frac{2M}{3M-8\pi}\geq 2/3$  for $M>8\pi$ and $f^{-1}$ is decreasing,
thus  the lower bound \eqref{tcbound2} on $T^*_c(n_0)$  follows immediately. The approximation  $f^{-1}(\frac{2}{3})\sim 0.874\,21$ 
is  obtained using any  simple numerical calculation software.
\vskip0.3cm

(2)\, 
The upper bound in \eqref{tcbound2bis} follows from \eqref{tcbound1} and the fact that
$$
\frac{1}{f^{-1}(a)}
=
\frac{1}{f^{-1}(\frac{2M}{3M-8\pi})}
\leq 
\frac{1}{\ln\left(1+\frac{M-8\pi}{2M}\right)}
=
\frac{1}{\ln\left(\frac{1}{a}  \right)},
$$
where 
$a=\frac{2M}{3M-8\pi}=\left( 1+\frac{M-8\pi}{2M}\right)^{-1}\in (0,1)$.
Indeed, if we set $f^{-1}(a)=\lambda>0$, i.e. $f(\lambda)=a$, then it is enough to show that
$$
\ln \left(\frac{\lambda}{1-e^{-\lambda}}\right)
\leq \lambda,\; \lambda>0.
$$
This is equivalent to the trivial inequality $1+\lambda\leq e^{\lambda}$ for $\lambda\geq 0$. So, the upper bound  \eqref{tcbound2bis} follows.
\vskip0.3cm

The lower bound of \eqref{tcbound2bis} is straightforward from \eqref{lowerbd2bis} of Theorem \ref{lowe1}
since  $M=  \pi\vert\vert n_0\vert\vert_{\infty} R^2$.
It can also be obtained  from \eqref{lowerbd2}
with  a well-chosen  finite $p$.
 Indeed, the function  $n_0$  is in all  $L^p$ and its $L^p$-norm  is explicitly computable.
More precisely, we have $ \vert \vert n_0\vert\vert_p=\sigma (\pi R^2)^{1/p}$ for all $1\leq p\leq +\infty$. So, we get 
$M:= \vert \vert n_0\vert\vert_1=\sigma (\pi R^2)$ and
$$
\left[ \frac{M}{  \vert \vert n_0\vert\vert_p}
\right]^{p'}
=\pi R^2, \; 1\leq p\leq +\infty,
$$
which proves  \eqref{tcbound2bis} by applying \eqref{lowerbd2}  with  $p$ chosen  large enough.
\vskip0.3cm

(3)
One can easily prove that $f^{-1}(v)\sim 2(1-v)$ as $v\rightarrow 1^-$.
Now since  we have $\lim_{M\rightarrow 8\pi^+} \frac{2M}{3M-8\pi}=1^-$, this implies that
$$
f^{-1}\left(  \frac{2M}{3M-8\pi}\right)\sim 
 \frac{M-8\pi}{8\pi}, \;\; {\rm as}\;\;  M\rightarrow 8\pi^+.
$$
Finally, we conclude the estimates \eqref{tcbound4}  and 
 \eqref{tcbound5} since $\pi\sigma R^2=M$ and  $M\rightarrow 8\pi^+$.
\vskip0.3cm

This concludes the proof of Proposition \ref{ball}. 
\hfill $\square$
\vskip0.3cm

Next, we compare the  estimates of $T^*_c(n_0)$  that can be obtained from the corollaries of Section \ref{tcsection} 
with the results of Proposition  \ref{ball}. 
\vskip0.3cm

$\bullet$ The first  alternative to obtain the upper bound of \eqref{tcbound2bis} is by applying  \eqref{tctrois} of Corollary  \ref{corapplica1}
(or equivalently in this case \eqref{tdeuxc}). 
This can be done because we can only consider the case $z_0=0$ by translation invariance of $T^*_c(n_0)$, 
and by the fact that $n_0$ is a non-increasing  radially symmetric function when $z_0=0$.
In our situation, we can explicitly compute $T^*_{c_2}(z)$ for  $z=0$. For this purpose, we compute 
$g_0^{\leftarrow}(t)$ for  $t\in [0,1)$.
We first easily  obtain
$$
g(\rho)
=\frac{{\mathcal M}_0(\rho)}
{M}
=\left[\inf (1,\rho/R)\right]^2, \, \rho>0.
$$
Hence, we get 
$g_0^{\leftarrow}(t)=R\sqrt{t}$ for all  $t\in (0,1)$.
By taking $t=a^{\theta}\in (0,1)$ where $a=\frac{2M}{3M-8\pi}$,
it follows that
$$
T^*_{c_2}(0)
=\frac{1}
{
4\ln\left(
1+\frac{M-8\pi}{2M}
\right)
}.
\inf_{\theta \in (0,1)}
\left(
\frac{ 
\left[
g_0^{\leftarrow}(a^{\theta})
\right]^2
}{1-\theta}
\right)
=
\frac{R^2}{4\ln\left(
1+\frac{M-8\pi}{2M}
\right)}
.
\inf_{\theta \in (0,1)}
\left(
\frac{ 
a^{\theta}
}{1-\theta}
\right).
$$
It is easily  shown that
$\inf_{\theta \in (0,1)}
\left(
\frac{ 
a^{\theta}
}{1-\theta}\right)=
\lim_{\theta \rightarrow 0^+} 
\left(
\frac{ 
a^{\theta}
}{1-\theta}\right)
=1$ using the fact that $a\in (2/3,1)$, hence $a\geq e^{-1}$.
Finally, we  conclude  that 
$$
T^*_c\leq 
T^*_{c_2}(0)=\frac{R^2}{4\ln\left(
1+\frac{M-8\pi}{2M}
\right)}.
$$
This proves once  again the upper bound as in \eqref{tcbound2bis}.
\vskip0.3cm

$\bullet$
A second  possibility of getting the exact value of $T^*_{c_2}(z_0)$ is to apply part  (1) of Corollary \ref{corapplica3}.
As already said, by   translation  invariance we  only need to  consider the case  $z_0=0$.
From the expression of $g_0^{\leftarrow}(t)=R\sqrt{t}$ described above, we deduce that $h(t):=(g_0^{\leftarrow})^2(t)=R^2 t$ 
for $t\in (0,1)$.  (Note that this formula   is also valid  for $t\in [0,1]$).
It is easy to check the sufficient conditions of  part  (1) of Corollary \ref{corapplica3}.
Indeed, the function $h$ is continuous and $h^{\prime}=R ^2>0$.
The associated function $F$ is given by $F(X)=X+1$  is a non-decreasing.
 We check  that $F(0^+)=1\geq \ln (\frac{1}{a})$ since  we always have 
$ \frac{1}{a}\leq  \frac{3}{2}\leq e$. Thus, the inequality \eqref{tcF1} of Corollary \ref{corapplica3}
asserts that
$$
T^*_{c_2}(0)=
\frac{h(1^{-})}
{4 \ln(1/a)}
=
\frac{R^2}
{4 \ln(1/a)},
$$
which  is the expected result.
\vskip0.3cm

$\bullet$ 
Again, a third  alternative  to obtain the upper bound of \eqref{tcbound2bis} is direct application of  \eqref{tctrois} of Corollary  \ref{corapplica1}
because $n_0$ has compact support.
To get the value of  $T^*_{c_3}$ in  \eqref{tctrois}, 
we need to compute   $R_0$ with  $K= \overline{B(0,R)}$.
 This can be done as follows. 
We first  prove that  $i_K(z)=\vert z\vert +R$  for all $z\in \R^2$, and 
 we deduce that $R_0=\inf_{z\in K}i_K(z)= i_K(0)=R$.
 We  apply \eqref{tctrois} of Corollary \ref{corapplica1}  to finally obtain 
$$
T^*\leq T_{c_3}^*:=
\frac{R^2}{4\ln\left(
1+\frac{M-8\pi}{2M}
\right)}.
$$
A more heuristic and geometric proof  to see that $R_0=R$ in \eqref{tctrois} is to find  
the smallest closed disk containing  $K= \overline{B(0,R)}$, the support of $n_0$,
which  is obviously $K=\overline{B(0,R)}$ itself.
\vskip0.3cm

\begin{rem}
The lower and upper bound of $T^*_c(n_0)$ in  \eqref{tcbound2bis} show that Corollary  \ref{corapplica1} and   (1) of Corollary   \ref{corapplica3},
and also Theorem \ref {lowe1} are sharp for estimating $T^*_c(n_0)$ for the class of  initial data $n_0$ of  characteristic functions of disks  
(up to universal multiplicative constants).
More precisely, each  result of Corollary  \ref{corapplica1} or Theorem \ref {lowe1}  shows the sharpness of each other. 
 A similar remark also holds with part  (1) of Corollary \ref{corapplica3} and  Theorem \ref {lowe1}.  
\end{rem}

\begin{rem}
From the upper bound of $T^*$ by $T_{c_3}^*$ described just above,  we  deduce   the next asymptotic estimate
$$
T^*\leq T_{c_3}^*
\sim
\frac{2M^2}
{4\pi\sigma (M-8\pi)}
\sim 
\frac{2(16\pi)}
{\sigma (M-8\pi)}, \quad M\rightarrow8\pi^+,
$$
 for fixed $\sigma>0$.
This asymptotic behaviour of $T_{c_3}^*$ can be  compared with 
the  asymptotic behaviour  of $T^*_c(n_0)$ obtained in \eqref{tcbound5}  as $M\rightarrow8\pi^+$ (for fixed $\sigma>0$).
In this case, it follows  immediately that $T_{c_3}^*\sim 2T^*_c(n_0)$. Thus, the upper estimate \eqref{tctrois} of $T^*$ in Corollary \ref{corapplica1}
is twice greater  than  the upper estimate of   $T^*$ obtained from $T^*_c(n_0)$   in \eqref{tcbound5}.
The same remark holds true  with  the estimate  \eqref{tcbound4}  instead of   \eqref{tcbound5}  when now $R>0$ is fixed and 
$M\rightarrow8\pi^+$. This is deduced  from the next asymptotic estimate
$$
T_{c_3}^*
\sim
\frac{4\pi R^2}
{M-8\pi}, \quad M\rightarrow8\pi^+,
$$
 for fixed $R>0$.
\end{rem}

\begin{rem}
For practical purposes, the quantity  $f^{-1}(\frac{2M}{3M-8\pi})$ can certainly be evaluated by numerical approximations  for any given mass $M>8\pi$.
So, the critical time  $T^*$ can  explicitly be  bounded  using  the inequality  \eqref{tcbound1}.
\end{rem}

We conclude this section with   an improvement for  the upper bound on $T^*_c(n_0)$
by applying the inequality  \eqref{tstarparti2} of Corollary  \ref{corapplica2}.
Indeed, it is easy to show that  the $2$-variance of the initial data $n_0$ is finite and 
 given by $V_{2}(n_0)={R^2}/{2}$.
Thus, we have  the next result.

\begin{pro}
Let $n_0=\sigma {\bf 1}_{B(z_0,R)}$,  for fixed $z_0\in \Ri^2$ and $\sigma, R>0$.
Assume that $M=\pi \sigma R^2>8\pi$. Then we have 
\begin{equation}\label{tcboundbetter}
T^ * \leq T^*_c(n_0) 
\leq
T^ *_{c_4}=
\frac{2^{-1} R^2}{4\ln\left(1+\frac{M-8\pi}{2M}\right)}.
\end{equation}
\end{pro}

\begin{rem}
The inequality  \eqref{tcboundbetter} clearly  improves by half
the upper bound of $ T^*_c(n_0) $  in the inequality \eqref{tcbound2bis} of Proposition \ref{ball}.
\end{rem}

\subsection{More  examples of initial data}
Many other examples could be presented  in this paper,
but we shall limit ourselves to a few of them in this last section.
\vskip0.3cm

\begin{enumerate}
\item
{\it Initial data uniformly supported by an annulus}

For instance,  we can consider the situation where  the cells are  uniformly distributed on an annulus 
in the plane corresponding to the following initial data,
 $$
n_0(x)=\sigma\left( \mathbbm{1}_{B(z_0,R_2)} (x)- \mathbbm{1}_{B(z_0,R_1)}(x) \right),  \, x\in \R^2,
$$
where $\sigma>0$ and $0<R_1<R_2$.
For this example,  we can prove similarly to the case of  the  characteristic function of a disk the following  result,
 \begin{equation}\label{hbound}
 T^* \leq   T^*_c(n_0) \leq 
 \frac{1}{4h^{-1}\left(\frac{L(M)}{\sigma\pi}\right)},
 \end{equation}
where
$h^{-1}$ the inverse of the  function 
$
h(s)=s^{-1}\left[ e^{-R_1^2s}-e^{-R_2^2s}\right], s>0,
$
 and 
$L(M)=\frac{2M^2}{3M-8\pi}$.
This example is an example of a $z_0$-radially symmetric initial data with compact support but not non-increasing. 
Similar results as in Proposition \ref{ball} obtained for  the case of a characteristic function of a disk could be given.
Let just mention at least two explicit results about $T^ *_c(n_0)$:
$$
e^{-1} 
\frac{R_2^2-R_1^2}
{4\ln\left(1+\frac{M-8\pi}{2M}\right)}
\leq T^ *_c(n_0)
\leq
\frac{R_2^2+R_1^2}{8\ln\left(1+\frac{M-8\pi}{2M}\right)}.
$$
The upper bound is obtained from an easy  computation of the variance $V_2(n_0)$.
We shall not provide the details here.
\vskip0.3cm

Other  examples  of   radially symmetric   (but not  necessarily non-increasing)  initial  data  $n_0$  
 for  which the inverse function of  the Laplace transform  is explicitly  computable
can be treated by applying Theorem \ref{tclaplace}.
 For instance, we can consider  the following limited  list  of examples.
 Recall that  $a>0$  is defined by the next  formula:
 $\frac{1}{a}=\frac{M}{L(M)}=1+\frac{M-8\pi}{2M}$.

\item
{\it Initial data  of the polynomial-Gaussian   form}.
Let $n_0$ defined by
$$
n_0(x)=
\sigma \vert x-z_0\vert^{2n}
e^{-\alpha \vert x-z_0\vert^{2}},
\quad x\in \R^2,
$$
where $n\in \N$, $\sigma, \alpha>0$ and $z_0\in \R^2$. Let $M= \sigma \pi n! \alpha^{-n-1}>8\pi$.
By applying Theorem \ref{tclaplace}, we can prove that
\begin{equation}\label{tcpolyn}
 T^*\leq T^*_c(n_0)
 \leq 
 (4\alpha)^{-1}  \left[ \left(\frac{1}{a}\right)^{\frac{1}{n+1}} -1 \right]^{-1}.
 \end{equation}

\item
{\it Initial data obtained by difference of two Gaussian functions}.
Let $z_0\in \R^2$, $0<d<b<+\infty$ and  $\sigma >0$ be fixed.
Let
$$
n_0(x)=\frac{\sigma}{b-d}\left( e^{-d\vert x-z_0\vert^2}- e^{-b\vert x -z_0\vert^2} \right), \; x\in \R^2,
$$
and let  $M=\frac{\pi \sigma}{db}>8\pi$  (i.e. $\sigma >8db$).
By applying Theorem \ref{tclaplace}, we can prove that
\begin{equation}\label{tcdoublexp}
T^*\leq T^*_c(n_0) \leq 
\left[
2\left(
\sqrt{ (b-d)^2+  \frac{4bd}{a}}
-
(b+d)
\right)
\right]^{-1}.
\end{equation}
As already seen for the other examples treated in this paper, the computation of the  2-variance is rather easy.
We have  $V_2(n_0)=\frac{b+d}{bd}$. 
By applying Corollary \ref{corapplica2}, we get immediately
\begin{equation}\label{tcdoublexpvar}
T^*\leq T^*_c(n_0) \leq 
T_{c_4}^*=
\frac{
\left(
\frac{b+d}{bd}
\right)
}
{4\ln\left(1+\frac{M-8\pi}{2M}\right)}.
\end{equation}
\end{enumerate}

We do not provide here the details of computations for the examples introduced in this  short section, 
nor a detailed comparison of the results that could  be deduced from 
Theorem \ref{lowe1} and  Corollaries   \ref{corapplica0}, \ref{corapplica1}, \ref{corapplica3}, \ref{corapplica2}   to avoid a lengthy paper.


\end{document}